\documentclass[12pt]{article}

\usepackage{graphics,amsmath,amssymb,amsthm,mathrsfs}

\newcommand{\rn}[1]{\mathbb{R}^{#1}}
\newcommand{\pdydx}[2]{\frac{\partial #1}{\partial #2}}

\newtheorem{thm}{Theorem}[section]
\newtheorem{lemma}[thm]{Lemma}

\theoremstyle{definition}
\newtheorem{remark}[thm]{Remark}

\numberwithin{equation}{section}

\begin{document}

\title{The $L^p$ Regularity Problem on Lipschitz Domains}

\author{Joel Kilty \and Zhongwei Shen}

\date{}

\maketitle

\begin{abstract}
\noindent This paper contains two results on the $L^p$ regularity problem
on Lipschitz domains. For second order elliptic systems and $1<p<\infty$,
we prove that the solvability
of the $L^p$ regularity problem is equivalent to that of the $L^{p^\prime}$
Dirichlet problem. For higher order elliptic equations and systems,
we show that if $p>2$, the solvability of the $L^p$ regularity problem
is equivalent to a weak reverse H\"older condition with exponent $p$.

\end{abstract}

\bibliographystyle{amsplain}

\section{Introduction}

Let $\Omega$ be a bounded Lipschitz domain in $\rn{d}$ with
connected boundary. Consider the elliptic system of
order $2\ell$, $\mathcal{L}(D)u=0$ in $\Omega$, where
$u=(u^1,\dots, u^m)$,
\begin{eqnarray}
    (\mathcal{L}(D)u)^j &=& \sum_{k=1}^m \mathcal{L}^{jk}(D)u^k,
    \mbox{\indent} j=1,\dots, m, \label{ellipticOperator1}\\
    \mathcal{L}^{jk}(D) &=& \sum_{|\alpha| = |\beta|=\ell}
    a^{jk}_{\alpha\beta} D^{\alpha}D^{\beta},
    \label{ellipticOperator2}
\end{eqnarray}
and $D=(D_1,D_2,\dots, D_d)$, $D_i=\partial / \partial
x_i$ for $i=1,2,\dots, d$.  Also, $\alpha =
(\alpha_1,\alpha_2,\dots, \alpha_d)$ is a multi-index with length
$|\alpha| = \alpha_1+\dots + \alpha_d$ and
$D^{\alpha}=D_1^{\alpha_1}D_2^{\alpha_2}\dots D_d^{\alpha_d}$.  Let

\begin{equation}
    \mathcal{L}^{jk}(\xi) =\sum_{|\alpha|=|\beta|=\ell}
    a^{jk}_{\alpha\beta} \xi^{\alpha} \xi^{\beta} \mbox{\indent for
    } \xi \in \rn{d}.
\end{equation}

\noindent We will assume throughout this paper that the
$a^{jk}_{\alpha\beta}$ are real constants satisfying the
Legendre-Hadamard ellipticity condition
\begin{equation} \label{ellipticityCondition}
    \mu |\xi|^{2\ell}|\eta|^2 \leq \sum_{j,k=1}^{m}
    \mathcal{L}^{jk}(\xi)\eta^j\eta^k \leq \frac{1}{\mu}
    |\xi|^{2\ell}|\eta|^2,
\end{equation}
for some $\mu>0,$ and all $\xi\in\rn{d}, \eta\in \rn{m}$,
as well as the symmetry condition
\begin{equation} \label{symmetryCondition}
    \mathcal{L}^{jk}(\xi) = \mathcal{L}^{kj}(\xi).
\end{equation}

The $L^p$ Dirichlet problem for the elliptic system
$\mathcal{L}(D)u=0$ in $\Omega$ consists of finding a solution
$u$ such that $(\nabla^{\ell-1} u)^*\in L^p(\partial\Omega)$ and $u, \nabla u, \dots, \nabla^{\ell-1} u$
take the prescribed data on $\partial\Omega$ in the sense of nontangential convergence.
Here and thereinafter $\nabla^k u$ denotes the tensor of all partial derivatives of order
$k$ and
$(w)^*$ the nontangential maximal function of $w$.
More precisely, let $\mbox{\emph{WA}}^{k,p} (\partial\Omega, \rn{m})$ denote the completion of
the set of arrays of functions
\begin{equation}\label{array}
\big\{ \dot{f}=(f_\alpha)_{|\alpha|\le k}
=(D^\alpha f)_{|\alpha|\le k}:\ f\in C_0^\infty(\rn{d}, \rn{m})\big\},
\end{equation}
under the scale-invariant norm on $\partial\Omega$,
\begin{equation} \label{pk norm}
\|\dot{f}\|_{k,p}:= \sum_{|\alpha|\le k} |\partial\Omega|^{\frac{k-|\alpha|}{1-d}}
\| D^\alpha f\|_p,
\end{equation}
where $\|\cdot\|_p$ denotes the norm in $L^p(\partial\Omega)$.
The $L^p$ Dirichlet problem is said to
be uniquely solvable if given any $\dot{f}\in
\mbox{\emph{WA}}^{\ell-1,p}(\partial\Omega, \rn{m})$,
there exists a unique function $u$ such that
\begin{equation} \label{DirichletProblem}
    \left\{\begin{array}{ll} \mathcal{L}(D)u =0 & \mbox{ in }
    \Omega,
    \\ D^{\alpha} u = f_{\alpha} & \mbox{ on } \partial\Omega
    \mbox{\ \ \ for } |\alpha|\leq \ell -1, \\
    (\nabla^{\ell-1} u)^* \in L^p(\partial
    \Omega). & \end{array} \right.
\end{equation}
Moreover, the solution $u$ satisfies the estimate
\begin{equation}\label{DirichletEstimate}
\| (\nabla^{\ell-1} u)^*\|_p \le C\, \sum_{|\alpha|=\ell-1} \| f_\alpha\|_p.
\end{equation}

If the Dirichlet data
in (\ref{DirichletProblem}) are taken from
$\mbox{\emph{WA}}^{\ell, p}(\partial\Omega, \rn{m})$ instead of \\$\mbox{\emph{WA}}^{\ell-1, p}(\partial\Omega, \rn{m})$, then they all have tangential
derivatives in $L^p(\partial\Omega)$.
Consequently  we may expect the solution to have
one order higher regularity. This is the so-called $L^p$ regularity problem.
Let $\nabla_t g$ denote the tangential derivatives of
$g$ on $\partial\Omega$.
We say that the $L^p$ regularity problem for $\mathcal{L}(D)u=0$ in
$\Omega$ is uniquely solvable if given $\dot{f}=\{ f_\alpha: |\alpha|\le \ell\}\in
{\mbox{\emph{WA}}}^{\ell,p}(\partial\Omega,\rn{m})$, there exists a unique function $u$
such that
\begin{equation} \label{RegProblem}
    \left\{\begin{array}{ll} \mathcal{L}(D)u =0 & \mbox{ in }
    \Omega,
    \\ D^{\alpha} u = f_{\alpha} & \mbox{ on } \partial\Omega
    \mbox{\ \ \ for } |\alpha|\leq \ell -1, \\
    (\nabla^{\ell} u)^* \in L^p(\partial
    \Omega). & \end{array} \right.
\end{equation}
Moreover, the solution $u$ satisfies the estimate
\begin{equation} \label{regEstimate}
    \|(\nabla^{\ell}u)^*\|_p \leq C \sum_{|\alpha|=\ell-1}
    \|\nabla_t f_{\alpha}\|_p.
\end{equation}

For $p$ close to $2$ and $d\ge 2$,
 the solvability of the $L^p$ Dirichlet and regularity problems
was established in \cite{dahlberg3,fabes2,fabes,gao} for second order elliptic systems and in \cite{dkv:biharmonic,pv:higherOrder,verchota:polyharmonic,verchota:higherOrder}
for higher order elliptic equations and systems. In the lower dimensional case $d=2$ or $3$,
the $L^p$ Dirichlet problem was solved for $2-\varepsilon<p\le \infty$ and
the $L^p$ regularity problem
for $1<p<2+\varepsilon$ (both ranges are sharp)
in \cite{dahlberg4,pv:biharmonic,pv:maximum,verchota:higherOrder}.
In the higher dimensional case $d\ge 4$,
the $L^p$ Dirichlet problem for $2<p< \frac{2(d-1)}{d-3} +\varepsilon$
was recently solved by Shen in \cite{shen:ellip,shen:necsuff}
for higher-order elliptic equations and systems. The paper \cite{shen:ellip} also established
the solvability of the $L^p$ regularity problem
for the second order elliptic systems in the case  $d\ge 4$ and $\frac{2(d-1)}{d+1} -\varepsilon <p<2$.
Related results may be found in \cite{kilty,shen:stokes,wright}
for the Stokes system and in \cite{shen:biharmonic}
for the biharmonic equation.
We remark that the results mentioned above extend the classical work
of Dahlberg, Jerison, Kenig, and Verchota in \cite{dahlberg2,dahlberg1,dahlberg,jerison,verchota} on $L^p$ boundary value problems
for Laplace's equation in Lipschitz domains.

In this paper we establish two related results on the $L^p$ regularity problem.
First, for general higher order elliptic equations and systems in $\Omega$,  we show
that if $p>2$, the solvability of the $L^p$ regualrity problem is equivalent
to a weak reverse H\"older condition with exponent $p$ on $\partial\Omega$.
Let $\Delta(P,r)= B(P,r)\cap \partial\Omega$ where $P\in\partial \Omega$.
The result may be formulated as follows.

\begin{thm} \label{necSuffThm}
    Let $\mathcal{L}(D)$ be a system of elliptic operators of order $2\ell$ satisfying conditions
    (\ref{ellipticityCondition}) and (\ref{symmetryCondition}).  For any bounded Lipschitz
    domain $\Omega$ and $p>2$, the following are equivalent.

    \begin{enumerate}
        \item The $L^p$ regularity problem for $\mathcal{L}(D)u=0$ in
 $\Omega$ is uniquely solvable.
        \item There exist $C>0$ and $r_0>0$ such that for any $P\in\partial\Omega$ and $0<r<r_0$, the
        weak reverse H\"{o}lder condition
        \begin{equation} \label{reverseHolder}
            \left(\frac{1}{r^{d-1}} \int_{\Delta(P,r)} |(\nabla^{\ell}
            v)^*|^p\,d\sigma \right)^{\frac{1}{p}} \leq C
            \left(\frac{1}{r^{d-1}} \int_{\Delta(P,2r)}
            |(\nabla^{\ell}v)^*|^2\,d\sigma\right)^{\frac{1}{2}},
        \end{equation}
 holds for any solution $v$ of $\mathcal{L}(D)v=0$ in $\Omega$
        with the properties $(\nabla^{\ell} v)^*\in L^2(\partial\Omega)$ and $D^{\alpha}v=0$
        on $\partial\Omega$ for $|\alpha|\leq \ell-1$ on $\Delta(P,3r)$.
    \end{enumerate}
\end{thm}

Theorem \ref{necSuffThm} extends the work of Shen in \cite{shen:necsuff} where a similar result was established
for the $L^p$ Dirichlet problem. When combined with the main result in
\cite{shen:necsuff}, it yields the following.

\begin{thm}\label{corollary}
    Let $\Omega$ be a bounded Lipschitz domain in $\rn{d}$, $d\geq
    4$.  Suppose that the $L^p$ regularity problem
for $\mathcal{L}(D)u=0$ in $\Omega$ is uniquely solvable
    for some $2<p<d-1$.
    Then the $L^q$ Dirichlet problem for $\mathcal{L}(D)u=0$
in $\Omega$ is uniquely solvable for $2<q<q_0+\varepsilon$, where
 $\frac{1}{q_0}=\frac{1}{p}-\frac{1}{d-1}$.
\end{thm}

In the second part of this paper we consider the case of second order
elliptic systems, i.e. $\ell=1$.  In this special case we show that
for any given Lipschitz domain,
the $L^p$ regularity problem and the $L^{p^\prime}$ Dirichlet problem are
in fact equivalent.

\begin{thm} \label{DirichletRegularityEquiv}
    Let $1<p<\infty$ and $\Omega$ be a bounded Lipschitz domain.
 Then, for any second order elliptic system
satisfying conditions (1.4)-(1.5),  the $L^{p}$ regularity
    problem in $\Omega$ is uniquely solvable if and only if the $L^{p^\prime}$ Dirichlet
    problem in $\Omega$ is uniquely solvable.
\end{thm}

We point out that although it is not implicitly stated, the duality between the
regularity and Dirichlet problems was essentially
established in the case of star-shaped Lipschitz domains
for Laplace's equation in \cite{verchota}.
Some partial results on this duality relation may also be found in \cite{shen:dirichletRegularity}
for second order elliptic equations with bounded measurable coefficients.
Our approach to the elliptic systems uses the basic duality argument in \cite{verchota}.
The main contribution here is a localization argument which allows us to treat the case of general
Lipschitz domains in the absence of positivity.

It would be very interesting to see if the duality between
the regularity and Dirichlet problems in Theorem 1.3 extends
to higher order elliptic equations and systems.
As a first step in this direction, some partial results have been obtained by the authors
for the biharmonic equation $\Delta^2 u=0$.

Note that by Theorems \ref{necSuffThm} and \ref{DirichletRegularityEquiv},
for second order elliptic systems, as in the case $p>2$,
the solvability of the $L^p$ Dirichlet problem
for $1<p<2$ is also equivalent to a weak reverse H\"older condition.
In particular, by the well known self-improving property of weak
reverse H\"older conditions,
it follows that if the Dirichlet problem in $\Omega$ is solvable
for some $1<p<2$, then it is also solvable in $\Omega$ for some $1<\bar{p}<p$.
This, together with Theorem \ref{corollary} as well as results in
\cite{dahlberg4}, gives the following.

\begin{thm} \label{corollary1} Consider the second order elliptic
system $\mathcal{L}(D)u=0$ in a bounded Lipschitz domain $\Omega$ in $\rn{d}$.
Let $\mathcal{A}$ denote the set of exponents $p\in (1,\infty)$
for which the Dirichlet problem in $\Omega$
is uniquely solvable.
Then

\begin{enumerate}

\item
If $d=2$ or $3$, $\mathcal{A}=(q, \infty)$ where $1\le q<2$.

\item

If $d\ge 4$, $\mathcal{A}=(q, s)$ where $1\le q<2<s\le \infty$.
Furthermore, $s=\infty$ if $q\le \frac{d-1}{d-2}$,
and $s\ge \frac{q(d-1)}{qd-2q -d +1}$ if $q>\frac{d-1}{d-2}$.

\end{enumerate}

\end{thm}

The paper is organized as follows. In Section 2 we show that
the condition $(\nabla^{2\ell-1} u)^*\in L^p$ for some $p>1$ implies that
$\nabla^{2\ell-1} u$ has nontangential limits on $\partial\Omega$ and
$\|(\nabla^{2\ell-1} u)^*\|_p \le C\, \| \nabla^{2\ell-1} u\|_p$.
The proof of Theorem \ref{necSuffThm} is given in Sections 3 and 4,
while Theorem \ref{corollary} is proved in Section 5. Sections 6 and 7 are
devoted to the proof of Theorem \ref{DirichletRegularityEquiv}.

Finally we remark that the summation convention will be used throughout this paper.
Also, $\Omega$ will always be a bounded Lipschitz domain with connected
boundary. We will use $\Gamma(x)=(\Gamma_{jk}(x))_{m\times m}$ to denote
a matrix of fundamental solutions on $\rn{d}$ to the operator $\mathcal{L}(D)$
with pole at the origin.

\section{A preliminary estimate}

\begin{thm} \label{nonTangLimitLemma}
   Suppose that $\mathcal{L}(D)u=0$ in $\Omega$ and
 $(\nabla^{2\ell -1} u)^*\in L^p(\partial\Omega)$ for
    some $p>1$. Then $\nabla^{2\ell-1} u$ has nontangential limits
    a.e. on $\partial\Omega$.  Furthermore,
    $\nabla^{2\ell-1} u \in L^p(\partial\Omega)$ and
\begin{equation}\label{nontangentialestimate}
\|(\nabla^{2\ell-1} u)^*\|_p \leq
    C\, \|\nabla ^{2\ell-1} u\|_p,
\end{equation}
where $C$ depends only on $\mu$, $d$, $m$, $\ell$,
 $p$ and the Lipschitz character of $\Omega$.
\end{thm}

\begin{proof}
    Let $\{ \Omega_r\}$ be a sequence of smooth domains such that
    $\Omega_r\uparrow \Omega$.
Fix $x\in \Omega$.  Since $\mathcal{L}(D)u=0$ in $\Omega$
    and $\mathcal{L}(D) \Gamma =\delta_0$, we may write
    \begin{eqnarray}
        D^{\gamma+ k } u^s(x) &=& \int_{\Omega_r}
        a_{\alpha\beta}^{ij} D^{\alpha +\beta} \Gamma^x_{js}(y)\cdot
        D^{\gamma +k}u^i\,dy \\ && \qquad\qquad+ \int_{\Omega_r}
        a^{ij}_{\alpha\beta} D^{\gamma +k} \Gamma^x_{js}(y)\cdot
        D^{\alpha +\beta} u^i\,dy, \nonumber
    \end{eqnarray}
where $\Gamma^x(y) =\Gamma(x-y)$, and
$\gamma$, $k$ are two
    multi-indicies with $|\gamma|=\ell$ and $|k|=\ell-1$.

    Next, we derive the Green's representation formula by
    integrating by parts to switch the derivatives on $\Gamma^x$ and
    $u$. This produces only boundary terms as the solid integrals
    cancel out.  Note that we should move derivatives in such a way
that no more than $2\ell$ derivatives are taken on either
    $\Gamma^x$ or $u$.  By doing so we obtain
    \begin{equation} \label{greenRepApprox}
        D^{\gamma+k} u^s(x) = \sum_{|\alpha|=2\ell-1}
        \int_{\partial\Omega_r} D^{\alpha} \Gamma^x_{js}\cdot \Pi_{ij}^\alpha
        (u^i)\,d\sigma,
    \end{equation}
 where $\Pi_{ij}^\alpha (u^i) $ is a sum of derivatives of $u^i$ of order $|\alpha|$
times various components of the unit normal to
$\partial\Omega_r$.

Let $\Lambda_r: \partial\Omega\to \partial\Omega_r$ denote the homeomorphism given by Theorem 1.12 in
\cite{verchota}.
We now rewrite (\ref{greenRepApprox}) as an
    integral on $\partial\Omega$ to obtain
    \begin{equation}
        D^{\gamma+k} u^s(x) = \sum_{|\alpha|=2\ell-1}
        \int_{\partial\Omega} D^{\alpha} \Gamma_{js}^x(\Lambda_r(P))\cdot
        \Pi_{ij}^\alpha (u^i)(\Lambda_r(P)) \omega_r\,d\sigma.
    \end{equation}
Since
 $$
|\Pi_{ij}^\alpha (u^i) (\Lambda_r (P))| \leq
    C\, (\nabla^{2\ell-1}u)^*(P),
$$
it follows
    that $$\|\Pi_{ij}^\alpha(u^i)(\Lambda_r)\|_p \leq
    C\, \|(\nabla^{2\ell-1}u)^*\|_p<\infty.$$
Consequently, there exists a
    subsequence, which we still denote by $\Pi_{ij}^\alpha (u^i) (\Lambda_r)$,
 that converges weakly in
    $L^p(\partial\Omega)$ to some function $g_j^\alpha\in L^p(\partial\Omega)$ as
$r\to\infty$.  Since
    $D^{\alpha} \Gamma^x_{js}(\Lambda_r(P)) \rightarrow D^{\alpha} \Gamma_{js}^x(P)$
    uniformly on $\partial\Omega$, we obtain
    $$
D^{\gamma+k}u^s(x) =
        \sum_{|\alpha|=2\ell-1} \int_{\partial\Omega} D^{\alpha}
        \Gamma_{js}^x(P)\cdot g_j^\alpha\,d\sigma.
$$
 This implies that $D^{\gamma +k}u$ has
    nontangential limits a.e. on $\partial\Omega$.
As a result, we have $| g_j^\alpha (P)|\le C|\nabla^{2\ell -1} u(P)|$
for a.e. $P\in \partial\Omega$.
It follows that
    $$\|(D^{\gamma +k}u^s)^*\|_p\le C\, \sum_{|\alpha|=2\ell-1} \sum_j
\| g_j^\alpha\|_p
\leq C\|\nabla^{2\ell-1} u\|_p.$$
This finishes the proof.
\end{proof}

\section{Sufficiency of the  weak reverse H\"{o}lder condition}

The goal of this section is to show that given any Lipschitz
domain $\Omega$ and any $p>2$,
the weak reverse H\"{o}lder condition (\ref{reverseHolder}) is
sufficient for the solvability of the $L^p$ regularity problem on
$\Omega$.

\begin{thm} \label{sufficiencyTheorem}
    Let $\mathcal{L}(D)$ be an elliptic operator of order $2\ell$
   satisfying the ellipticity condition
    (\ref{ellipticityCondition}) and the symmetry condition
    (\ref{symmetryCondition}).  Let $\Omega$ be a bounded Lipschitz
    domain in $\rn{d}$ and fix $p>2$.  Suppose that for any
    $\Delta(P,r)\subset \partial\Omega$ with $P\in \partial\Omega$ and
$0<r<r_0$, the reverse H\"{o}lder
    condition (\ref{reverseHolder}) holds for all solutions of $\mathcal{L}(D)u=0$ in
    $\Omega$ with the properties $(\nabla^{\ell} u)^*\in
    L^2(\partial\Omega)$ and $D^{\alpha}u=0$ on $\Delta(P,3r)$ for
    $|\alpha|\leq \ell -1$.  Then the $L^p$ regularity problem in $\Omega$ is
    uniquely solvable.
\end{thm}

 The proof of the following Poincar\'{e}
type inequality may be found in
\cite{shen:necsuff}.

\begin{lemma} \label{shenPoincare}
    Suppose $\ell \geq 1$.  Let $\dot{f}=\{f_{\alpha}: |\alpha|\leq
    \ell\} \in {\mbox{\emph{WA}}}^{2,\ell}(\partial\Omega)$ and $\Delta(P,r)\subset
    \partial\Omega$.  Then, there exists a polynomial $h$ of degree
    at most $\ell -1$ such that
    \begin{equation} \label{shenPoincareEstimate}
        \left\| f_{\beta} - D^{\beta}h\right\|_{L^2(\Delta(P,r))}
        \leq Cr^{\ell-|\beta|}\sum_{|\alpha|=\ell-1} \|\nabla_t
        f_{\alpha} \|_{L^2(\Delta(P,r))},
    \end{equation}
 for any multi-index $\beta$ with $|\beta|\leq \ell
    -1$.
\end{lemma}

The proof of Theorem \ref{sufficiencyTheorem}, which is similar to that of
Theorem 3.1 in \cite{shen:necsuff},  relies on a real
variable argument. The proof of the following theorem may be found in \cite{shen:boundary}.

\begin{thm}\label{realVariableArgument}
    Let $S=\left\{ (x',\psi(x')):x'\in \rn{d-1}\right\}$ be a Lipschitz
    graph in $\rn{d}$. Let $Q_0$ be a surface cube in  $S$ and $F\in L^2(2Q_0)$.  Let
    $p>2$ and $g\in L^q(2Q_0)$ for some $2<q<p$.  Suppose that for
    each dyadic subcube $Q$ of $Q_0$ with $|Q|\leq \beta |Q_0|$,
    there exists two integrable functions $F_Q$ and $R_Q$ on $2Q$
    such that $|F|\leq |F_Q|+|R_Q|$ on $2Q$, and
    \begin{eqnarray}
        \lefteqn{\left( \frac{1}{|2Q|} \int_{2Q} |R_Q|^p\,d\sigma
        \right)^{1/p} \leq } \nonumber\\&& \qquad C_1\left\{\left(\frac{1}{|\alpha Q|} \int_{\alpha
        Q} |F|^2\,d\sigma\right)^{1/2} + \sup_{Q'\supset Q} \left(\frac{1}{|Q'|}
        \int_{Q'} |g|^2\,d\sigma\right)^{1/2}\right\}, \label{realVariableAssump1}\\
        \lefteqn{\left(\frac{1}{|2Q|}\int_{2Q} |F_Q|^2\,d\sigma\right)^{1/2} \leq
        C_2\sup_{Q'\supset Q} \left(\frac{1}{|Q'|} \int_{Q'}
        |g|^2\,d\sigma\right)^{1/2},} \label{realVariableAssump2}
    \end{eqnarray}
 where $C_1,C_2>0$ and $0<\beta<1<\alpha$.  Then,
    \begin{eqnarray} \label{shenThmConclusion}
        \left( \frac{1}{|Q_0|}\int_{Q_0} |F|^q\,d\sigma\right)^{1/q}
        &\leq& C_3\left\{\left(\frac{1}{|2Q_0|}\int_{2Q_0} |F|^2\,d\sigma\right)^{1/2}\right. \\ && \left. \qquad\qquad+
        \left(\frac{1}{|2Q_0|}\int_{2Q_0}
        |g|^q\,d\sigma\right)^{1/q} \right\}, \nonumber
    \end{eqnarray}
 where $C_3$ depends only on $d$, $p$, $q$, $C_1$,
    $C_2$, $\alpha$, $\beta$, and $\|\nabla \psi\|_{\infty}$.
\end{thm}

We now proceed to the proof of Theorem \ref{sufficiencyTheorem}.

\noindent{\bf Proof of Theorem \ref{sufficiencyTheorem}.}

    The uniqueness for $p>2$ follows from the uniqueness for $p=2$.
    To establish the existence, we let $\dot{f}=\{ f_\alpha: |\alpha|\le \ell\}
\in {\mbox{\emph{WA}}}^{\ell,p}(\partial\Omega)$ and $u$ be
    the solution to the $L^2$ regularity problem with boundary data
    $\{ f_\alpha: |\alpha|\le \ell -1\}$.  We  will show that if $P\in\partial\Omega$ and
    $0<s\le cr_0$, then
    \begin{eqnarray}
        \left(\frac{1}{s^{d-1}}\int_{\Delta(P,s)}
        |(\nabla^{\ell}u)^*|^p\,d\sigma \right)^{\frac{1}{p}}\nonumber
        &\leq&  C
        \left(\frac{1}{s^{d-1}} \int_{\Delta(P,Cs)}
        |(\nabla^{\ell}u)^*|^2\,d\sigma\right)^{\frac{1}{2}}  \\ && \qquad+
        C\left(\frac{1}{s^{d-1}} \int_{\Delta(P,Cs)}
        \sum_{|\alpha|=\ell-1} |\nabla_t f_{\alpha}|^p\,d\sigma
        \right)^{\frac{1}{p}}. \label{reverseHolderClaim}
    \end{eqnarray}
 By covering $\partial\Omega$ with a finite number of
    balls of radius $cr_0$, estimate (\ref{reverseHolderClaim})
    implies that
    \begin{eqnarray*}
        \|(\nabla^{\ell} u)^*\|_p &\leq&
        C|\partial\Omega|^{\frac{1}{p}-\frac{1}{2}} \|(\nabla^{\ell}
        u)^*\|_2 + C \sum_{|\alpha|=\ell -1} \|\nabla_t
        f_{\alpha}\|_p \\
        &\leq& C|\partial\Omega|^{\frac{1}{p}-\frac{1}{2}}
        \sum_{|\alpha|=\ell -1} \|\nabla_t f_{\alpha}\|_2 +
        C\sum_{|\alpha|=\ell -1} \|\nabla_t f_{\alpha}\|_p \\
        &\leq& C\sum_{|\alpha|=\ell -1} \|\nabla_t f_{\alpha}\|_p.
    \end{eqnarray*}
Here we have used the $L^2$
    regularity estimate as well as H\"{o}lder's inequality.  We now
    seek to establish estimate $(\ref{reverseHolderClaim})$.

    Fix $P\in \partial\Omega$ and $0<s<cr_0$.  By rotation and
    translation we may assume that $P=0$ and
    \begin{eqnarray}
        B(P,r_0) \cap \Omega &=& B(P,r_0) \cap \big\{ (x',x_d)\in
        \rn{d}: x_d >\psi(x') \big\}, \label{Lipschitz-graph-1}\\
        B(P,r_0) \cap \partial\Omega &=& B(P,r_0) \cap \big\{
        (x',\psi(x'))\in\rn{d}: x'\in \rn{d-1} \big\},\label{Lipschitz-graph-2}
    \end{eqnarray}
    \noindent where $\psi:\rn{d-1}\rightarrow \rn{}$ is a Lipschitz
    function.  Consider the surface cube
$$Q_0 =\big\{ (x',\psi(x'))\in \rn{d-1}\times
    \rn{}: \ |x_1|<s,\dots, |x_{d-1}|<s\big\}.$$
 Let $Q$ be a small subcube of $Q_0$ with diameter $r$.
   Choose $\varphi \in C^{\infty}_0(\rn{d})$ with $0\leq \varphi
    \leq 1$, $\varphi=1$ on $8Q$, $\varphi=0$ on
    $\rn{d}\backslash 16Q$, and $|D^{\alpha}\varphi| \leq
    \frac{C}{r^{|\alpha|}}$ for $|\alpha| \leq 2\ell$.  Let $h$ be
    the polynomial of degree at most $\ell-1$, given by
 Lemma \ref{shenPoincare}, but with $\Delta(P,r)$
    replaced with $16Q$.  Write $u=v+w+h$ where $v$ is the
    solution to the $L^2$ regularity problem with boundary data
    $(u-h)\varphi$ and $w$ is the solution to the $L^2$ regularity
    problem with boundary data $(1-\varphi)(u-h)$.  Note that for
    $|\alpha| \leq \ell-1$
    \begin{equation}
        D^{\alpha} v = D^{\alpha}((u-h)\varphi) = \sum_{\beta\leq
        \alpha} \frac{\alpha !}{\beta!(\alpha-\beta)!}
        (f_{\beta}-D^{\beta}h)D^{\alpha-\beta}\varphi.
    \end{equation}
Now, let
    \begin{eqnarray*}
        F &=& |(\nabla^{\ell} u)^*|, \\
        g &=& \sum_{|\alpha|=\ell -1} |\nabla_t f_{\alpha}|, \\
        F_Q &=& 2|(\nabla^{\ell}v)^*|, \\
        R_Q &=& 2|(\nabla^{\ell}w)^*|.
    \end{eqnarray*}
 Using the fact that $v$ is the solution of the
    $L^2$ regularity problem, we obtain
    \begin{eqnarray}
        \frac{1}{|2Q|} \int_{2Q} |F_Q|^2\,d\sigma
        &\leq&  \frac{C}{|2Q|} \int_{\partial\Omega}
        |(\nabla^{\ell} v)^*|^2\,d\sigma \nonumber \\
        &\leq& \frac{C}{|2Q|}\sum_{|\alpha|=\ell-1}
        \int_{\partial\Omega} |\nabla_t D^{\alpha} v|^2\,d\sigma.
        \label{tanGradEst1}
    \end{eqnarray}
    \noindent Now, note that
    \allowdisplaybreaks{
    \begin{eqnarray}
        & \ &\sum_{|\alpha|=\ell -1} \int_{\partial\Omega} |\nabla_t
        D^{\alpha} v|^2\,d\sigma \nonumber\\
 &\leq&
        C\sum_{|\beta|\le \ell-1}\int_{16Q} |\nabla_t (f_{\beta}-D^{\beta}
        h)|^2\frac{d\sigma}{r^{2(\ell-1-|\beta|)}}
        \nonumber \\ && \qquad\qquad+C \sum_{|\beta|\le \ell-1}
        \int_{16Q} |f_{\beta}-D^{\beta}h|^2\frac{d\sigma}{r^{2(\ell -
        |\beta|)}} \nonumber \\
        &\leq& C\sum_{|\alpha|=\ell-1} \int_{16Q} |\nabla_t
        f_{\alpha}|^2\,d\sigma \nonumber \\ && \qquad\qquad + C \sum_i \sum_{|\beta|\leq \ell-2}
        \int_{16Q}
        |D_{x_i}(f_{\beta}-D^{\beta}h)|^2\frac{d\sigma}{r^{2(\ell-1-|\beta|)}}
        \nonumber \\ && \qquad\qquad + C \sum_{|\beta|\leq \ell-1}
        \frac{1}{r^{2(\ell-|\beta|)}}\int_{16Q}
        |f_{\beta}-D^{\beta}h|^2\,d\sigma \nonumber\\
        &\leq& C\sum_{|\alpha|=\ell-1} \int_{16Q} |\nabla_t
        f_{\alpha}|^2\,d\sigma + C \sum_{|\beta|\leq \ell-1}
        \frac{1}{r^{2(\ell-|\beta|)}}\int_{16Q}
        |f_{\beta}-D^{\beta}h|^2\,d\sigma \nonumber \\
        &\leq& C\sum_{|\alpha|=\ell -1} \int_{16Q} |\nabla_t
        f_{\alpha}|^2\,d\sigma, \label{tanGradEst2}
    \end{eqnarray}
    }
    \noindent where we have used Lemma \ref{shenPoincare} in the last
    step.  By combining estimates (\ref{tanGradEst1}) and
    (\ref{tanGradEst2}) we obtain
    \begin{equation}
        \frac{1}{|2Q|} \int_{2Q} |F_Q|^2\,d\sigma \leq
        \frac{C}{|16Q|}\sum_{|\alpha|=\ell-1} \int_{16Q}
        |\nabla_t f_{\alpha}|^2\,d\sigma \leq \frac{C}{|16Q|}
        \int_{16Q} |g|^2\,d\sigma.
    \end{equation}
    \noindent This implies that
    \begin{equation} \label{thmCond1}
        \left(\frac{1}{|2Q|} \int_{2Q} |F_Q|^2\,d\sigma
        \right)^{\frac{1}{2}} \leq C \sup_{Q'\supset Q}
        \left(\frac{1}{|Q'|} \int_{Q'}
        |g|^2\,d\sigma\right)^{\frac{1}{2}}.
    \end{equation}

    Note that $w$ is a solution of the $L^2$ regularity problem with
    $(\nabla^{\ell} w)^*\in L^2(\partial\Omega)$ and $D^{\alpha}w=0$
    on $16Q$ for $|\alpha|\leq \ell -1$.  Thus, we may use the
    weak reverse H\"{o}lder inequality (\ref{reverseHolder}) and the
    above estimates on $v$ to obtain
    \allowdisplaybreaks{
    \begin{eqnarray}
      &\ &  \left(\frac{1}{|2Q|} \int_{2Q}
        |R_Q|^p\,d\sigma\right)^{\frac{1}{p}}\nonumber =
        2\left(\frac{1}{|2Q|} \int_{2Q} |(\nabla^{\ell}
        w)^*|^p\,d\sigma\right)^{\frac{1}{p}} \nonumber \\
        &\ & \le C\left(\frac{1}{|4Q|} \int_{4Q}
        |(\nabla^{\ell} w)^*|^2\,d\sigma\right)^{\frac{1}{2}}
        \nonumber \\
        &\ & \le C\left(\frac{1}{|4Q|}\int_{4Q} |(\nabla^{\ell}
        u)^*|^2\,d\sigma\right)^{\frac{1}{2}} +
        C\left(\frac{1}{|4Q|} \int_{4Q} |(\nabla^{\ell}
        v)^*|^2\,d\sigma\right)^{\frac{1}{2}} \nonumber \\
        &\ & \leq C\left(\frac{1}{|4Q|}\int_{4Q}
        |F|^2\,d\sigma\right)^{\frac{1}{2}} + C\sup_{Q'\supset Q}
        \left(\frac{1}{|Q'|}\int_{Q'}
        |g|^2\,d\sigma\right)^{\frac{1}{2}}. \label{thmCond2}
    \end{eqnarray}
    }

    We should point out that the weak reverse H\"{o}lder condition
    on surface balls is equivalent to the weak reverse H\"{o}lder condition
    on surface cubes.  This is because we may cover a surface cube
    by sufficiently small surface balls with finite overlap and
    vice versa.  Thus, both conditions of Theorem
    \ref{realVariableArgument} are satisfied and estimate
    (\ref{reverseHolderClaim}) follows by covering $\Delta(P,s)$
    with a finite number of sufficiently small surface cubes.
    This establishes the solvability of the $L^q$ regularity problem
for any $2<q<p$. Finally, since the weak reverse H\"older condition
(\ref{reverseHolder}) has the self-improving property,
the argument above also gives the solvability of the $L^q$ regularity
problem for $2<q<p+\varepsilon$ and in particular, for $q=p$.
\qed

\section{Necessity of the weak reverse H\"{o}lder condition}

In this section we show that the reverse H\"{o}lder condition
(\ref{reverseHolder}) with exponent $p>2$ is also necessary for the
solvability of the $L^p$ regularity problem.

\begin{thm} \label{necessityTheorem}
    Let $\mathcal{L}(D)$ be an elliptic operator of order $2\ell$
    satisfying the ellipticity condition
    (\ref{ellipticityCondition}) and the symmetry condition
    (\ref{symmetryCondition}).  Let $\Omega$ be a bounded Lipschitz
    domain in $\rn{d}$ and fix $p>2$.  Suppose that the $L^p$
    regularity problem for $\mathcal{L}(D)u=0$ in $\Omega$ is uniquely
    solvable.  Then the weak reverse H\"{o}lder inequality
    (\ref{reverseHolder}) holds for solutions of $\mathcal{L}(D)u=0$
    in $\Omega$ with the properties $(\nabla^{\ell}u)^*\in
    L^2(\partial\Omega)$ and $D^{\alpha} u=0$ on $\Delta(P,3r)$ for
    $|\alpha|\leq \ell-1$.
\end{thm}

\begin{proof}
We begin by choosing $r_0>0$ so that for any $P\in \partial\Omega$,
(\ref{Lipschitz-graph-1})-(\ref{Lipschitz-graph-2}) hold
 after a possible rotation of the coordinate system.
    Fix $P_0 \in \partial\Omega$ and $0<r<cr_0$.  Let $u$ be a
    solution of $\mathcal{L}(D)u=0$ in $\Omega$ such that
    $(\nabla^{\ell}u)^*\in L^2(\partial\Omega)$ and $D^{\alpha}u=0$
    on $\Delta(P_0,10r)$ for $|\alpha|\leq \ell -1$.  For a function
    $v$ on $\Omega$ and $P\in \partial\Omega$ define
    \begin{eqnarray*}
        \mathcal{M}_1(v)(P) &=& \sup \{ |v(x)|: x\in \gamma(P)
        \mbox{ and } |x-P| < c_0r \}, \\
        \mathcal{M}_2(v)(P) &=& \sup \{ |v(x)|: x\in \gamma(P)
        \mbox{ and } |x-P| \geq c_0 r \},
    \end{eqnarray*}
where $$
\gamma(P)=\gamma_a(P)
=\{ x\in \Omega: \ |x-P|<(1+a) \text{dist}(x, \partial\Omega)\}
$$
and  $a>1$ is sufficiently large.
    \noindent Then, $(\nabla^{\ell} u)^* = \max \{
    \mathcal{M}_1(\nabla^{\ell}u), \mathcal{M}_2(\nabla^{\ell}u)
    \}$.  If $x\in \gamma(P)$ for some $P\in \Delta(P_0,r)$ and
    $|x-P|\geq c_0 r$, then by the interior estimates we have
$$
    |\nabla^{\ell} u(x)| \leq \frac{C}{r^{d}} \int_{B(x,cr)}
    |\nabla^{\ell}
    u(y)|\,dy \leq \frac{C}{r^{d-1}} \int_{\Delta(P_0,2r)}
    |(\nabla^{\ell} u)^*|\,d\sigma.
$$
It follows that for any $p>2$,
    \begin{equation} \label{vEstFar}
         \left(\frac{1}{r^{d-1}} \int_{\Delta(P_0,r)}
         |\mathcal{M}_2(\nabla^{\ell} u)|^{p}\,d\sigma
         \right)^{\frac{1}{p}} \leq C \left(\frac{1}{r^{d-1}}
         \int_{\Delta(P_0,2r)} |(\nabla^{\ell}
         u)^*|^2\,d\sigma\right)^{\frac{1}{2}}.
    \end{equation}

    We now estimate $\mathcal{M}_1(\nabla^{\ell} u)$ on
    $\Delta(P_0,r)$.  First, choose
    $\varphi \in C^{\infty}_0(\rn{d})$ such that $\varphi=1$ on
    $B(P_0,2r)$, $\varphi=0$ on $\rn{d}\backslash B(P_0,3r)$, and
    $|D^{\alpha} \varphi| \leq \frac{C}{r^{|\alpha|}}$ for $|\alpha|\leq
    2\ell$.  Note that
    \allowdisplaybreaks{
    \begin{eqnarray}
        \lefteqn{[ \mathcal{L}(D)(u\varphi)]^j  \sum_{k=1}^{m}
        \sum_{|\alpha|=|\beta|=\ell} a_{\alpha\beta}^{jk}
        D^{\alpha}D^{\beta} (u^k\varphi)} \nonumber \\
        \ \ \ \ &=& \sum_{k=1}^{m} \sum_{|\alpha|=|\beta|=\ell}
        a^{jk}_{\alpha\beta} D^{\alpha}\left\{ \varphi D^{\beta}u^k
        + \sum_{\gamma\leq \beta} \frac{\beta!}{\gamma!
        (\beta-\gamma)!} D^{\gamma} u^k D^{\beta-\gamma} \varphi
        \right\} \nonumber \\
        &=& \sum_{k=1}^{m} \sum_{|\alpha|=|\beta|=\ell}
        \sum_{\gamma<\alpha} a_{\alpha\beta}^{jk}
        \frac{\alpha!}{\gamma!(\alpha-\gamma)!} D^{\beta+\gamma}u^k
        D^{\alpha-\gamma} \varphi \label{expandOper} \\ && \qquad+ \sum_{k=1}^m
        \sum_{|\alpha|=|\beta|=\ell} \sum_{\gamma<\beta}
        a_{\alpha\beta}^{jk} \frac{\beta!}{\gamma!(\beta-\gamma)!}
        \sum_{\delta \leq \alpha}
        \frac{\alpha!}{\delta!(\alpha-\delta)!}
        D^{\beta-\gamma+\alpha-\delta}\varphi D^{\delta+\gamma} u^k,
        \nonumber
    \end{eqnarray}
    }
 where we used the fact that $\mathcal{L}(D)u=0$ in
    $\Omega$.

    Recall that $\Gamma(x)=(\Gamma_{ij}(x))_{m\times m}$ denotes a matrix of
    fundamental solutions on $\rn{d}$ to the operator
    $\mathcal{L}(D)$ with pole at the origin.  We remark that if $d$ is
    odd or $2\ell<d$, $\Gamma_{ij}(x)$ is homogeneous of degree $2\ell -
    d$ and smooth away from the origin.  If $d$ is even and $2\ell
    >d$, then $\Gamma_{ij}(x)=\Gamma_{ij}^{(1)}(x)+\ln{|x|} \cdot \Gamma_{ij}^{(2)}(x)$
    where $\Gamma_{ij}^{(1)}(x)$ is homogeneous of degree $2\ell-d$ and
    $\Gamma_{ij}^{(2)}(x)$ is a polynomial of degree $2\ell-d$.  In this case we
    replace $\ln{|x|}$ with $\ln{(|x|/r)}$.  This can be done since
    $\Gamma_{ij}^{(2)}(x)$ is a polynomial of degree $2\ell-d$.  In
    either case we have the estimate
    \begin{equation} \label{fundSolnDeriv}
        |D^{\alpha} \Gamma (x)| \leq
        \frac{C_{\alpha}}{|x|^{d-2\ell+|\alpha|}} \mbox{ \indent for
        } |\alpha| \geq 2\ell-d+1,
    \end{equation}
    \noindent since the derivatives $D^{\alpha}$ eliminate the logarithmic
    singularity if $|\alpha|> 2\ell-d$.

    Fix $y_0 \in \rn{d}\backslash \overline{\Omega}$ so that
    $|y_0-P_0|=r \sim \mbox{dist}(y_0,\partial\Omega)$.
As in \cite{shen:necsuff}, let
    $\Gamma(x,y)=\Gamma (x-y)$ and define
    \begin{equation} \label{adjFundSoln}
        F_{ij}(x,y) = \Gamma (x,y)- \sum_{|\gamma| \leq 2\ell-1}
        \frac{(y-y_0)^{\gamma}}{\gamma!} D_y^{\gamma}
        \Gamma_{ij}(x,y_0).
    \end{equation}
    \noindent The summation term in (\ref{adjFundSoln}) is a solution of
    $\mathcal{L}(D)u=0$ in $\Omega$ in both $x$ and $y$
    variables.  It is subtracted from $\Gamma(x,y)$ to create the
    desired decay when $|x-P_0|\geq 5r$ and $|y-P_0|\leq 4r$.
Let $T(P,s)=\Omega\cap B(P,s)$. By
    the Taylor remainder theorem and (\ref{fundSolnDeriv}), if $x\in
    \Omega \backslash T(P_0,5r)$ and $y\in T(P_0,3r)$, then
    \begin{equation} \label{adjFundSolnDerivAway}
        |\nabla^{\ell}_x D^{\alpha}_y F_{ij}(x,y)| \leq
        \frac{Cr^{2\ell - |\alpha|}}{|x-y|^{d+\ell}} \mbox{
        \indent for } |\alpha| \leq 2\ell .
    \end{equation}
    \noindent Also, if $x\in T(P_0,5r)$ and $y\in T(P_0,3r)$ we have
    \begin{equation} \label{adjFundSolnDerivNear}
        |\nabla^{\ell}_xD^{\alpha}_y F_{ij}(x,y)| \leq
        \frac{Cr^{\ell - |\alpha|-1}}{|x-y|^{d-1}} \mbox{ \indent
        for } |\alpha| \leq \ell -1.
    \end{equation}
    \noindent To establish (\ref{adjFundSolnDerivNear}) we consider
    two cases: $|\alpha|> \ell -d$ and $|\alpha| \leq \ell -d$.  For
    the first case we use estimate (\ref{fundSolnDeriv}).  The
    second case is handled by noting that the term involving the
    possible logarithmic function in $\nabla_x^{\ell}D^{\alpha}_y
    \Gamma_{ij}(x,y)$ is bounded by $C|x-y|^{\ell -d -
    |\alpha|}\ln \left| \frac{|x-y|}{r}\right|$.  Since $|x-y|\leq
    Cr$ and $\ell-|\alpha|-1>0$, it is bounded by the right hand side of
    (\ref{adjFundSolnDerivNear}).

    Next, define $w(x)=(w^1(x),\dots, w^m(x))$ by
    \begin{eqnarray*}
        w^i(x) &=&\sum_{|\alpha|=|\beta|=\ell}
        \sum_{\gamma<\alpha} (-1)^{\gamma}a^{jk}_{\alpha\beta}
        \frac{\alpha!}{\gamma!(\alpha-\gamma)!} \int_{\Omega}
        D_y^{\gamma}\left(F_{ij}(x,y)D^{\alpha-\gamma}\varphi \right)
        D^{\beta} u^k\,dy  \\   & & \quad\quad+
       \sum_{|\alpha|=|\beta|=\ell}
        \sum_{\gamma<\beta} \sum_{\delta \leq \alpha}
        (-1)^{\gamma} a_{\alpha\beta}^{jk}
        \frac{\beta!}{\gamma!(\beta-\gamma)!}\frac{\alpha!}{\delta!(\alpha-\delta)!}
        \times\\ && \quad \quad \quad \quad \quad \quad \quad \qquad\int_{\Omega}
        D_y^{\gamma}\left(F_{ij}(x,y)D^{\alpha-\delta+\beta-\gamma}\varphi\right)D^{\delta}u^k\,dy.
    \end{eqnarray*}
    Note that $\mathcal{L}(D)(w)=\mathcal{L}(D)(u\varphi)$ in
    $\Omega$.  This follows from integration by parts and
    (\ref{expandOper}).

    Now, on $\Delta(P_0,r)$ we have
$$
\mathcal{M}_1(\nabla^{\ell}u)
    = \mathcal{M}_1(\nabla^{\ell}(u\varphi)) \leq
    \mathcal{M}_1(\nabla^{\ell}w) +
    \mathcal{M}_1(\nabla^{\ell}(u\varphi -w)).
$$
For $x\in T(P_0,5r)$, we use (\ref{adjFundSolnDerivNear}) to obtain
    \begin{eqnarray} \label{localWEst}
        |\nabla^{\ell} w(x)| \leq  C \sum_{|\gamma|\le \ell}
       \frac{1}{r^{\ell-|\gamma|+1}} \int_{T(P_0,3r)\backslash T(P_0,2r)}
        \frac{|D^{\gamma}u(y)|}{|x-y|^{d-1}} \,dy.
    \end{eqnarray}
    Thus, if $P\in \Delta(P_0,r)$, we have
    \begin{eqnarray*}
        \mathcal{M}_1(\nabla^{\ell} w)(P) &\leq&
        C\sum_{|\gamma|\leq \ell}r^{|\gamma|-\ell-d}
        \int_{T(P_0,3r)}
        |D^{\gamma} u|\,dy
    \end{eqnarray*}
    
    \begin{eqnarray*}
        &\leq& C\left(\frac{1}{r^d}\int_{T(P_0,3r)}
        |\nabla^{\ell} u|^2\,dy\right)^{\frac{1}{2}} \\ && \qquad\qquad+
        C\sum_{|\gamma|\le \ell-1} r^{|\gamma|-\ell}
\left(\frac{1}{r^d} \int_{T(P_0,3r)}
        |D^{\gamma}u|^2\,dy\right)^{\frac{1}{2}} \\
&\leq& C\left(\frac{1}{r^d}\int_{T(P_0,3r)}
        |\nabla^{\ell} u|^2\,dy\right)^{\frac{1}{2}}\\
        &\leq& C\left(\frac{1}{r^{d-1}} \int_{\Delta(P_0,3r)}
        |(\nabla^{\ell}u)^*|^2\,d\sigma \right)^{\frac{1}{2}},
    \end{eqnarray*}
 where we have used the assumption $D^\gamma u =0$ on $\Delta (P_0, 10r)$ for  $|\gamma|\le \ell-1$ and
the Poincar\'e inequality in the third inequality.  This gives
    \begin{equation} \label{nonTangWEst}
        \left(\frac{1}{r^{d-1}}\int_{\Delta(P_0,r)}
        |\mathcal{M}_1(\nabla^{\ell}
        w)|^p\,d\sigma\right)^{\frac{1}{p}} \leq
        C\left(\frac{1}{r^{d-1}} \int_{\Delta(P_0,3r)}
        |(\nabla^{\ell} u)^*|^2\,d\sigma\right)^{\frac{1}{2}}.
    \end{equation}

    We still need to estimate
    $\mathcal{M}_1(\nabla^{\ell}(u\varphi-w))$.  This is where the
    assumption that the $L^p$ regularity problem on $\Omega$ is
    uniquely solvable is used.  Recall that $\mathcal{L}(D)(u\varphi-w)=0$
    in $\Omega$.  As in \cite{shen:necsuff},  we also have
    $(\nabla^{\ell-1}(u\varphi-w))^*\in L^2(\partial\Omega)$.
Thus, we may apply the uniqueness of the $L^2$ Dirichlet problem
    and the $L^p$ regularity estimate to obtain
    \begin{eqnarray*}
        \int_{\Delta(P_0,r)} |\mathcal{M}_1(\nabla^{\ell} (u\varphi
        -w))|^p\,d\sigma &\leq& \int_{\partial\Omega}
        |(\nabla^{\ell}(u\varphi-w))^*|^p\,d\sigma \\
        &\leq& C\int_{\partial\Omega} |\nabla_t\nabla^{\ell-1}
        (u\varphi -w)|^p\,d\sigma \\
        &\leq& C\int_{\partial\Omega} |\nabla_t
        \nabla^{\ell-1}(u\varphi)|^p\,d\sigma \\ && \qquad\qquad+
        C\int_{\partial\Omega} |\nabla_t\nabla^{\ell-1}
        w|^p\,d\sigma \\
        &=& C\int_{\partial\Omega} |\nabla_t \nabla^{\ell-1}
        w|^p\,d\sigma,
    \end{eqnarray*}
where the last step follows from the observation that $\nabla^{\ell-1} (u\varphi)=0$ on
    $\partial\Omega$.

    Now, let $p=\frac{q(d-1)}{d-q}$.  Then,
    $\frac{d-1}{p}=\frac{d}{q}-1$.  Using (\ref{localWEst}) and Lemma 4.2 in \cite{shen:necsuff} we obtain
    \begin{eqnarray*}
        {\left(\frac{1}{r^{d-1}}\int_{\Delta(P_0,5r)} |\nabla^{\ell}
        w|^p\,d\sigma\right)^{\frac{1}{p}}}
&\leq &
\frac{C}{r^{\frac{d-1}{p}}}
\left\| r^{|\gamma|-\ell -1}
\sum_{|\gamma|\le \ell} |D^\gamma u|\right\|_{L^q(T(P_0,3r))}\\
        &\leq& C \frac{1}{r^{\frac{d}{q}}}
        \left(\int_{T(P_0,3r)} |\nabla^\ell
        u|^q\,dy\right)^{\frac{1}{q}} \\
        &\leq& \left(\frac{1}{r^{d-1}} \int_{\Delta(P_0,3r)}
        |(\nabla^{\ell}u)^*|^{\bar{q}}\,d\sigma\right)^{\frac{1}{\bar{q}}},
    \end{eqnarray*}
 where $\bar{q}=\max \{ q,2\}$.  This gives
    \begin{equation} \label{wEstNear}
        \left(\frac{1}{r^{d-1}} \int_{\Delta(P_0,5r)} |\nabla^{\ell}
        w|^p\,d\sigma\right)^{\frac{1}{p}} \leq C
        \left(\frac{1}{r^{d-1}} \int_{\Delta(P_0,3r)}
        |(\nabla^{\ell} u)^*|^{\bar{q}}\,d\sigma
        \right)^{\frac{1}{\bar{q}}}.
    \end{equation}

    Finally, if $P\in \partial\Omega\backslash \Delta(P_0,5r)$,
we use estimate (\ref{adjFundSolnDerivAway}) to obtain
    \allowdisplaybreaks{
    \begin{eqnarray*}
        |\nabla^{\ell} w(P)| &\leq&
\frac{Cr^{d+\ell-1}}{|P-P_0|^{d+\ell -1}}
\sum_{|\gamma|\le \ell}
r^{|\gamma|-\ell-d}\int_{T(P_0,3r)} |D^\gamma u|\, dy\\
                &\leq& C
\frac{r^{d+\ell-1}}{|P-P_0|^{d+\ell -1}}
\left( \frac{1}{r^d} \int_{T(P_0,3r)} |\nabla^{\ell}
        u|^2\,dy\right)^{\frac{1}{2}} \\
        &\leq& \frac{Cr^{d+\ell-1}}{|P-P_0|^{d+\ell -1}} \left(\frac{1}{r^{d-1}}\int_{\Delta(P_0,3r)}
        |(\nabla^{\ell} u)^*|^2\,d\sigma\right)^{\frac{1}{2}}.
    \end{eqnarray*}
    }
    This implies that
    \begin{equation} \label{wEstFar}
        \left(\frac{1}{r^{d-1}} \int_{\partial\Omega\backslash
        \Delta(P_0,5r)} |\nabla^{\ell} w|^p\,d\sigma
        \right)^{\frac{1}{p}} \leq C\left(\frac{1}{r^{d-1}}
        \int_{\Delta(P_0,3r)} |(\nabla^{\ell}
        u)^*|^{\bar{q}}\,d\sigma\right)^{\frac{1}{\bar{q}}}.
    \end{equation}
Combining estimates (\ref{vEstFar}), (\ref{nonTangWEst}),(\ref{wEstNear}), and (\ref{wEstFar}),
    we have proved that
    \begin{equation} \label{qBarVEst}
        \left(\frac{1}{r^{d-1}} \int_{\Delta(P_0,r)} |(\nabla^{\ell}
        u)^*|^p\,d\sigma\right)^{\frac{1}{p}} \leq C
        \left(\frac{1}{r^{d-1}} \int_{\Delta(P_0,3r)}
        |(\nabla^{\ell}
        u)^*|^{\bar{q}}\,d\sigma\right)^{\frac{1}{\bar{q}}},
    \end{equation}
 where $\bar{q}=\max (q,2)$ and
    $\frac{d-1}{p}=\frac{d}{q}-1$. Note that since $p\ge 2$,
    $$
\frac{1}{q}-\frac{1}{p} =
    \frac{1}{d}\left(1-\frac{1}{p}\right) \geq \frac{1}{2d}.
$$
Thus    we may iterate estimate (\ref{qBarVEst}) to obtain
    \begin{equation}
        \left(\frac{1}{r^{d-1}} \int_{\Delta(P_0,cr)}
        |(\nabla^{\ell} u)^*|^p\,d\sigma \right)^{\frac{1}{p}}
        \leq C\left(\frac{1}{r^{d-1}} \int_{\Delta(P_0,r)}
        |(\nabla^{\ell} u)^*|^2\right)^{\frac{1}{2}},
    \end{equation}
starting with $q=2$.  We may do this since the $L^p$
    solvability of the regularity problem implies the $L^s$
    solvability for $2<s<p$.

    By covering $\Delta(P_0,r)$ with sufficiently small surface balls
    $\{ \Delta(P_j,cr)\}$, we obtain the weak reverse H\"{o}lder condition
    (\ref{reverseHolder}).
\end{proof}

\begin{remark}
Under the same assumptions as in Theorem \ref{necessityTheorem}, we have
    \begin{eqnarray}
        \lefteqn{\hspace{-1.0in}\left(\frac{1}{r^{d-1}} \int_{\Delta(P_0,r)}
        |\mathcal{M}_1(\nabla^{\ell}u)|^p\,d\sigma
        \right)^{\frac{1}{p}}} \nonumber \\ \hspace{0.5in}&\leq&
        \frac{C}{r}\left(\frac{1}{r^{d-1}}
        \int_{\Delta(P_0,6r)} |(\nabla^{\ell-1}
        u)^*|^2\,d\sigma\right)^{\frac{1}{2}},\label{remark-1}
    \end{eqnarray}
where $\mathcal{L}(D)u=0$ in $T(P_0,10r)$,
    $\mathcal{M}_1(\nabla^{\ell-1}u)\in L^2(\Delta(P_0,10r))$, and
    $D^{\alpha} u=0$ on $\Delta(P_0,10r)$ for $|\alpha|\leq \ell -1$.
Indeed, a careful inspection of the proof of Theorem \ref{necessityTheorem} shows that
    \begin{equation}\label{remark-2}
        \left(\frac{1}{r^{d-1}} \int_{\Delta(P_0,r)}
        |\mathcal{M}_1(\nabla^{\ell}u)|^p\,d\sigma\right)^{\frac{1}{p}}
        \leq C \left(\frac{1}{r^d}\int_{T(P_0, 3r)}
        |\nabla^\ell u|^2\,d\sigma \right)^{\frac{1}{2}}.
    \end{equation}
Using $D^\alpha u=0$ on $\Delta (P_0,10r)$ for $|\alpha|\le \ell-1$ and
the fact that the $L^2$ regularity problem is uniquely solvable on every
bounded Lipschitz domain, one may deduce that the right hand side of (\ref{remark-2})
is bounded by the right hand side of (\ref{remark-1}). We leave the details to the reader.
\end{remark}

\section{Proof of Theorem 1.2}

Let $\frac{1}{q_0}=\frac{1}{p}-\frac{1}{d-1}$. By Theorem 1.1 in \cite{shen:necsuff},
  it suffices to establish the weak reverse H\"{o}lder condition
    \begin{eqnarray}
        \lefteqn{\hspace{-0.5in} \left( \frac{1}{r^{d-1}} \int_{\Delta(P_0,r)}|(\nabla^{\ell-1} u)^*|^{q_0}\,d\sigma
        \right)^{\frac{1}{q_0}}} \nonumber\\ && \leq C \left(\frac{1}{r^{d-1}}
        \int_{\Delta(P_0,20r)} |(\nabla^{\ell-1}
        u)^*|^2\,d\sigma\right)^{\frac{1}{2}},\label{5-1}
    \end{eqnarray}
where $\mathcal{L}(D)u=0$ in $\Omega$, $(\nabla^{\ell-1} u)^* \in L^2(\partial\Omega)$ and
$D^\alpha u =0$ on $\Delta(P_0,100r)$ for $|\alpha|\le \ell-1$.
Clearly we may assume that $\Omega\cap B(P_0, Cr)$ is given by the region above a Lipschitz graph
in the sense of (\ref{Lipschitz-graph-1})-(\ref{Lipschitz-graph-2}).

    Let $P\in \Delta(P_0,r)$ and $x\in \gamma (P)$. It follows from the interior estimates
that
    \begin{equation}\label{5-2}
        |\nabla^{\ell}u(x)| \leq \frac{C}{s^{d-1}} \int_{|y-P|< cs}
        |(\nabla^{\ell}u)^*(y)|\,d\sigma,
    \end{equation}
\noindent where $s=\mbox{dist}(x,\partial\Omega)$.  Next, write $$D^{\alpha}u(x',x_d) -
    D^{\alpha}u(x',\tilde{x}_d) = - \int_{x_d}^{\tilde{x}_d}
    D^{\alpha+e_d}u(x',s)\,ds,$$  where $|\alpha|=\ell-1$.  This, together with (\ref{5-2}),
 gives that
    \begin{equation}
        \mathcal{M}_1(\nabla^{\ell-1}u)(P) \leq
        C\int_{\Delta(P_0,3r)}\frac{\widetilde{\mathcal{M}}_1 (\nabla^{\ell}
        u)}{|P-y|^{d-2}}\,d\sigma(y) +
        \widetilde{\mathcal{M}}_2(\nabla^{\ell-1}u) (P),
    \end{equation}
 where $\widetilde{\mathcal{M}}_1$ and $\widetilde{\mathcal{M}}_2$
 are defined in the same fashion as $\mathcal{M}_1$ and $\mathcal{M}_2$,
but using a family of slightly larger nontangential approach regions
$\{ \gamma_b(P): P\in \partial\Omega\}$, where $b>a$. Thus,
by the fractional integral estimates as well as the obvious estimate for $\widetilde{\mathcal{M}}_2
(\nabla^{\ell-1}u)$, we have
    \begin{eqnarray}
        \lefteqn{\hspace{-0.5in}\left(\frac{1}{r^{d-1}} \int_{\Delta(P_0,r)}
        |(\nabla^{\ell-1}u)^*|^{q_0}\,d\sigma\right)^{\frac{1}{q_0}}}
        \nonumber\\
        &\leq&
C\left(\frac{1}{r^{d-1}}\int_{\Delta(P_0, 3r)} |(\nabla^{\ell-1} u)^*|^2d\sigma\right)^{1/2}
\label{5-4} \\ && \qquad\qquad+ Cr
        \left(\frac{1}{r^{d-1}} \int_{\Delta(P_0,3r)}
        |\widetilde{\mathcal{M}}_1
(\nabla^{\ell}u)|^p\,d\sigma\right)^{\frac{1}{p}}. \nonumber
    \end{eqnarray}
The desired estimate (\ref{5-1}) now follows from (\ref{5-4}) and (\ref{remark-1}).

\section{Duality between the regularity and Dirichlet problems, part I}

The remaining two sections of this paper are devoted to the proof Theorem \ref{DirichletRegularityEquiv}.
In this section we show that for any second order elliptic system,
the solvability of the $L^p$ regularity problem implies
that of the $L^{p^\prime}$ Dirichlet problem.

To simplify the notation in the case $\ell=1$, we write the $m\times m$
system as
$\mathcal{L}(u)=0$, where
 $(\mathcal{L}(u))^\alpha =a_{ij}^{\alpha\beta}D_iD_j u^\beta$ for $\alpha=1,\dots, m,$
and
\begin{equation}\label{ellipticity}
\mu |\xi|^2 |\eta|^2\le a_{ij}^{\alpha\beta} \xi_i\xi_j \eta^\alpha \eta^\beta
 \le \frac{1}{\mu} |\xi|^2 |\eta|^2
\end{equation}
for some $\mu>0$, and all $\xi\in \rn{d}$ and $\eta\in \rn{m}$.
Without the loss of generality,
we may assume that $a_{ij}^{\alpha\beta}=a_{ji}^{\beta\alpha}$ in the place of (1.4).

\begin{thm} \label{regImpliesDirichlet}
    Let $1<p<\infty$ and $\Omega$ be a bounded Lipschitz domain.
   Suppose that the $L^{p}$ regularity
    problem for $\mathcal{L}(u)=0$ in $\Omega$ is uniquely solvable.
Then, the $L^{p^\prime}$ Dirichlet
    problem for $\mathcal{L}(u)=0$ in $\Omega$ is uniquely solvable.
\end{thm}

\begin{proof}
    We begin with the existence.
    Let $f\in C^{\infty}_0(\rn{d})$ and $u$ be the
    solution of the $L^{p}$ regularity problem in $\Omega$ with boundary
    data $f$; that is,
$
 \mathcal{L}(u)=0
$
in
$\Omega$, $u=f$ on $\partial\Omega$
and $\|(\nabla u)^*\|_p \le C\, \| \nabla_t f\|_p$.
    Since $(\nabla u)^*\in L^{p}(\partial\Omega)$, it follows from
Theorem \ref{nonTangLimitLemma} that
    $\nabla u$ exists a.e. on $\partial\Omega$ in the sense of nontangential convergence.
  Therefore, we have
    the Green's representation formula
    \begin{equation}
        u(x) = \int_{\partial\Omega}
        \Gamma^x\cdot\pdydx{u}{\nu}\,d\sigma - \int_{\partial\Omega}
        \pdydx{\Gamma^x}{\nu} \cdot u\,d\sigma,
    \end{equation}
where $\Gamma^x (y)=\Gamma(x-y)$ and
$\frac{\partial u}{\partial \nu}$ denotes the conormal derivative
of $u$ on $\partial\Omega$, defined by
\begin{equation}\label{conormal}
\left(\frac{\partial u}{\partial\nu}\right)^\alpha
=a_{ij}^{\alpha\beta}n_iD_j u^\beta,
\end{equation}
where $n=(n_1, \dots, n_d)$ is the outward unit normal to $\partial\Omega$.
Thus, by the well known singular integral estimates on Lipschitz surfaces,
    \begin{equation} \label{dirichletEstRegProb}
        \|(u)^*\|_{p^\prime} \leq C \left\|
        \pdydx{u}{\nu}\right\|_{-1,p^\prime} + C \|f\|_{p^\prime},
    \end{equation}
where $\| \cdot\|_{-1, p^\prime}$ denotes the norm in $W^{-1,p^\prime}(\partial\Omega)$, the
dual of the Sobolev space $W^{1,p}(\partial\Omega)$ equiped with the scale-invariant norm
$$
\|h\|_{W^{1,p}(\partial\Omega)}
=\|\nabla_t h\|_p +|\partial\Omega|^{\frac{1}{1-d}} \| h\|_p.
$$

    To estimate the first term on the right hand side of
    (\ref{dirichletEstRegProb}), we let $g\in
    C^{\infty}_0(\rn{d})$ and $w$ be the solution of the $L^{p}$
    regularity problem with data $g$.  Using integration by parts, we
    have
    \begin{eqnarray*}
        \left| \int_{\partial\Omega} \pdydx{u}{\nu}\cdot g \,d\sigma
        \right| = \left| \int_{\partial\Omega} u\cdot
        \pdydx{w}{\nu}\,d\sigma \right|
        \leq \|f\|_{p^\prime}\left\|\pdydx{w}{\nu}\right\|_{p}.
    \end{eqnarray*}
    \noindent Now, since $$ \left\|\pdydx{w}{\nu} \right\|_{p} \leq C
    \|(\nabla w)^*\|_{p} \leq C \|\nabla_t g\|_{p} \leq
    C\|g\|_{W^{1,p}(\partial\Omega)},$$
by duality we obtain
    \begin{equation} \label{dualityEst}
        \left\|\pdydx{u}{\nu}\right\|_{-1,p^\prime} \leq C \|f\|_{p^{\prime}}.
    \end{equation}
    \noindent So, by combining estimates (\ref{dirichletEstRegProb})
    and (\ref{dualityEst}), we obtain
   $
        \|(u)^*\|_{p^\prime} \leq C \|f\|_{p^\prime}.
    $

    Next, we establish the existence in the general case.  Let $f\in
    L^{p^\prime}(\partial\Omega)$ and choose $f_k\in C^{\infty}_0(\rn{d})$ such that
    $f_k\rightarrow f$ in $L^{p^\prime}(\partial\Omega)$ as $k\to \infty$.  Let $u_k$ be the unique solution
of the $L^p$ regularity problem in $\Omega$ with data $f_k$. Since $u_j-u_k$ is the unique
solution of the $L^p$ regularity problem  with data $f_j-f_k$,
 we have
    $$
\|(u_j-u_k)^*\|_{p^\prime}\leq C\| f_j-f_k\|_{p^\prime}\rightarrow 0.$$
    Using $$\sup_K |u_j-u_k| \leq C_K
    \|(u_j-u_k)^*\|_{p^\prime},
$$
 where $K\subset\subset \Omega$ is compact,
    we see that $u_j\rightarrow u$ uniformly on every compact subset of
    $\Omega$.  This implies that $\mathcal{L}(u)=0$ and
    $\|(u)^*\|_{p^\prime}\leq C \|f\|_{p^\prime}$.  To see this, note that
    $\|(u_j)^*_{\varepsilon}\|_{p^\prime} \leq C\|f_j\|_{p^\prime}$ where
    $$
(u_j)_\varepsilon^*(P)=\sup\{ |u_j(x)|: x\in \gamma(P) \text{ and dist}(x,\partial\Omega)
\geq\varepsilon\}.
$$
  Letting
    $j\rightarrow \infty$ and then $\varepsilon\rightarrow 0$ we
    obtain $\|(u)^*\|_{p^\prime}\leq C\|f\|_{p^\prime}$, as desired.

    To complete the existence part, we need to show that $u\rightarrow f$ nontangentially
    almost everywhere.  First, note that
    $\|(u_j-u_k)^*_{\varepsilon}\|_{p^\prime} \leq C\|f_j-f_k\|_{p^\prime}$.
    Letting $k\rightarrow \infty$ and then $\varepsilon\rightarrow
    0$, we obtain $$\|(u_j-u)^*\|_{p^{\prime}} \leq C \|f_j-f\|_{p^{\prime}}.$$  To show that
 $u$ has nontangential limits, we define
    $$\Lambda(P) = \limsup_{\substack{x\rightarrow P \\ x\in
    \gamma(P)}} u(x) -\liminf_{\substack{x\rightarrow P \\ x\in
    \gamma(P)}} u(x).$$
Now, fix $j$ and note that
    \begin{eqnarray*}
        \Lambda &=& \limsup (u-u_j+u_j) - \liminf (u-u_j+u_j) \\
        &\leq& \limsup (u-u_j) + \limsup u_j - \liminf (u-u_j)
        -\liminf u_j \\
        &=& \limsup (u-u_j) -\liminf (u-u_j) \\
        &\leq& 2(u-u_j)^*.
    \end{eqnarray*}
    Thus, $\|\Lambda\|_{p^\prime}\leq C\|(u-u_j)^*\|_{p^\prime}\rightarrow
    0$ as $j\rightarrow \infty$.  This implies that $\Lambda=0$
    a.e. and hence $u$ has nontangential limits a.e. on $\partial\Omega$.

    It remains to show that $u\rightarrow f$ a.e on
    $\partial\Omega$.  To do this, let
$$\widetilde{\Lambda}(P) =|
    \lim_{\substack{x\rightarrow P \\ x\in \gamma(P)}} u(x) -
    f(P)|.$$
Note that for any $j$,
    \begin{eqnarray*}
        \widetilde{\Lambda} = | \lim_{\substack{x\rightarrow P \\
        x\in \gamma(P)}} (u-u_j) + (f_j-f)(P)|
        \leq |(u-u_j)^*(P)| + |(f_j-f)(P)|.
    \end{eqnarray*}
This implies that
$$
\|\widetilde{\Lambda}\|_{p^\prime} \leq
    \|(u-u_j)^*\|_{p^\prime} + \|f_j-f\|_{p^\prime} \rightarrow 0
\quad \text{ as } j\to \infty.
$$
 Thus, $\widetilde{\Lambda}=0$ a.e. and $u=f$ a.e. on $\partial\Omega$ in the sense of nontangential
convergence.

    In the second part of this proof we establish the uniqueness of the solution.
  To this end, suppose that
$$\left\{ \begin{array}{ll}
    \mathcal{L}(u)=0 & \mbox{ in } \Omega, \\ u=0 & \mbox{ a.e. on }
    \partial\Omega, \\ (u)^*\in L^{p^\prime}(\partial\Omega). & \end{array} \right.$$
    We need to show that $u\equiv 0$ in $\Omega$.

    Fix $x\in \Omega$ and let $G^x=G(x,y)=\Gamma(x-y)-v^x(y)$, where
$v^x$ is the solution of the $L^p$ regularity problem with data $\Gamma(x-y)$; i.e.,
    $$
\left\{\begin{array}{ll}\mathcal{L}(v^x)=0 & \mbox{ in } \Omega, \\
    v^x(y)=\Gamma(x-y) & \mbox{ on } \partial\Omega, \\ (\nabla v^x)^*\in
    L^{p}(\partial\Omega). & \end{array}\right.
$$
Note that $(\nabla v^x)^*
    \in L^{p}(\partial\Omega)$ implies that $(\nabla_yG^x)^{*,\varepsilon}
    \in L^{p}(\partial\Omega)$ if $2\varepsilon< \text{dist}(x, \partial\Omega)$,
 where we have used the notation
$$
(w)^{*,\varepsilon}(P)
=\sup \big\{ |w(x)|: \ x\in \gamma(P) \text{ and dist}(x,\partial\Omega)<\varepsilon\big\}.
$$

    Choose $\varphi_{\varepsilon}\in C^{\infty}_0(\rn{d})$ such that $\varphi_{\varepsilon}
    =1$ on $\{ y\in \Omega: \mbox{dist}(y,\partial\Omega)\geq \varepsilon\}$,
    $\varphi_{\varepsilon}=0$ on $\Omega \backslash \{y\in \Omega: \mbox{dist}(y,\partial
    \Omega)\geq \varepsilon/2\}$, $|\nabla \varphi_{\varepsilon}|\leq \frac{C}{\varepsilon}$, and
    $|\nabla^2 \varphi_{\varepsilon}|\leq \frac{C}{\varepsilon^2}$. Note that
    \begin{equation}
      u(x)=  u\varphi_{\varepsilon}(x) = \int_{\Omega} G(x,y)\mathcal{L}(u
        \varphi_{\varepsilon})\,dy
    \end{equation}
and
    \begin{eqnarray*}
        (\mathcal{L}(u\varphi_{\varepsilon}))^{\alpha} &=& a_{ij}^{\alpha\beta} D_iD_j\{ u^{\beta}
        \varphi_{\varepsilon}\} \\
        &=& a_{ij}^{\alpha\beta} \left\{ D_iD_j u^{\beta}\cdot \varphi_{\varepsilon} + D_ju^{\beta}
        \cdot D_i \varphi_{\varepsilon} + D_iu^{\beta}\cdot D_j \varphi_{\varepsilon} + u^{\beta}\cdot
        D_iD_j\varphi_{\varepsilon} \right\} \\
        &=& a_{ij}^{\alpha\beta} D_ju^{\beta}\cdot D_i \varphi_{\varepsilon} + a_{ij}^{\alpha\beta}
        D_iu^{\beta}\cdot D_j\varphi_{\varepsilon} + a_{ij}^{\alpha\beta} u^{\beta}
        D_iD_j\varphi_{\varepsilon},
    \end{eqnarray*}

    \noindent where we used the fact that $\mathcal{L}(u)=0$ in $\Omega$.  This implies that
    \begin{eqnarray}
        u^{\alpha} (x) &=& \int_{\Omega} G_{\alpha\gamma}(x,y)a_{ij}^{\gamma\beta}
        \left\{ D_ju^{\beta}\cdot D_i\varphi_{\varepsilon} + D_iu^{\beta}\cdot D_j\varphi_{\varepsilon}
        + u^{\beta}\cdot D_iD_j\varphi_{\varepsilon}\right\}dy
\nonumber\\
        &=& -\int_{\Omega} \left\{ D_jG_{\alpha\gamma}(x,y)a_{ij}^{\gamma\beta} u^{\beta}D_i\varphi_{\varepsilon}
        +G_{\alpha\gamma}(x,y)a_{ij}^{\gamma\beta} u^{\beta}D_iD_j\varphi_{\varepsilon}\right.
\label{green}\\ && \qquad\qquad
        +\left.
        D_iG_{\alpha\gamma}(x,y)a_{ij}^{\gamma\beta} u^{\beta} D_j\varphi_{\varepsilon}\right\}
        \,dy.\nonumber
    \end{eqnarray}
It follows that
    \begin{eqnarray}
        |u^{\alpha}(x)| &\leq& \frac{C}{\varepsilon} \int_{E_\varepsilon}
        |\nabla_y G(x,y)||u(y)|\,dy + \frac{C}{\varepsilon^2} \int_{E_\varepsilon}
        |G(x,y)||u(y)|\,dy \nonumber\\
        &\leq& \frac{C}{\varepsilon} \left\{ \int_{E_\varepsilon} |\nabla_y G(x,y)|^{p}
        dy \right\}^{\frac{1}{p}} \left\{ \int_{E_\varepsilon} |u|^{p^\prime}\,dy \right\}^\frac{1}{p^\prime}
\nonumber \\ && \qquad+
        \frac{C}{\varepsilon^2}\left\{\int_{E_\varepsilon} |G(x,y)|^{p}\,dy\right\}^{\frac{1}{p}}
        \left\{\int_{E_\varepsilon} |u|^{p^\prime}\,dy \right\}^{\frac{1}{p^\prime}}, \label{uniqueEstGreen}
    \end{eqnarray}
where $E_\varepsilon=\{ x\in \Omega:(\varepsilon/2) \le \text{dist}(x, \partial\Omega)
\le \varepsilon\}$.

   Using the mean value theorem and $G(x, y)=0$ for $y\in \partial\Omega$,
it is easy to see that
 $|G(x,y)|\le C\varepsilon (\nabla_y G^x)^{*,\varepsilon} (P)$
if $y\in E_\varepsilon\cap \gamma (P)$.
It then follows from
    (\ref{uniqueEstGreen}) that
    \begin{eqnarray*}
        |u^{\alpha}(x)| &\leq& \frac{C}{\varepsilon^{1-\frac{1}{p}}} \left\{ \int_{\partial\Omega} |(\nabla_y
        G^x)^{*,\varepsilon}|^{p}\,d\sigma\right\}^{\frac{1}{p}}\left\{\int_{E_\varepsilon} |u|^{p^\prime}\,dy
        \right\}^{\frac{1}{p^\prime}} \\
        &\leq& C_x \left\{ \frac{1}{\varepsilon}\int_{E_\varepsilon} |u|^{p^\prime}\,
dy\right\}^{\frac{1}{p^\prime}} \\
        &\leq& C_x \left\{ \int_{\partial\Omega}
        |(u)^{*,\varepsilon}|^{p^\prime}\,d\sigma\right\}^{\frac{1}{p^\prime}}.
    \end{eqnarray*}

Finally, since $u=0$ a.e. on $\partial\Omega$ in the sense of nontangential convergence,
it follows that $(u)^{*,\varepsilon}\rightarrow 0$ a.e. as $\varepsilon \rightarrow 0$.
Using $(u)^{*,\varepsilon}\leq (u)^*\in L^{p^\prime}(\partial\Omega)$ and the dominated convergence theorem,
 we obtain
$\| (u)^{*,\varepsilon}\|_{p^\prime}\to 0
$
as $\varepsilon\rightarrow 0$.
 This implies that $u\equiv 0$ in $\Omega$ and the
    uniqueness of the solution is established.
\end{proof}

\section{Duality between the regularity and Dirichlet problems, part II}

    In this final section we prove the other implication in Theorem
    \ref{DirichletRegularityEquiv}.

\begin{thm} \label{DirichletImpliesRegularity}
    Let $\Omega$ be a bounded Lipschitz domain and $1<p<\infty$.  If
    the $L^{p^\prime}$ Dirichlet problem for $\mathcal{L}(u)=0$
in $\Omega$ is uniquely solvable, then the $L^{p}$
    regularity problem for $\mathcal{L}(u)=0$ in $\Omega$ is uniquely solvable.
\end{thm}

The proof of this theorem relies on the following lemma.
It also uses the solvability of the $L^2$ regularity problem for second order elliptic
systems, established in \cite{dahlberg3,fabes2,fabes,gao}.

\begin{lemma} \label{mainLemma}
Under the same assumptions as in Theorem \ref{DirichletImpliesRegularity}, we have
    \begin{equation} \label{mainestimate}
        \|(\nabla u)^*\|_p +|\partial\Omega|^{\frac{1}{1-d}} \| (u)^*\|_p
 \leq C\left\{ \|\nabla_t u\|_{p} +
        |\partial\Omega|^{\frac{1}{1-d}}\|u\|_{p} \right\},
    \end{equation}
where $u$ is the solution of the $L^2$ regularity problem with data
$f\in C^{\infty}_0(\rn{d})$.
\end{lemma}

We first demonstrate how to deduce Theorem \ref{DirichletImpliesRegularity} from  Lemma \ref{mainLemma}.

\noindent{\bf Proof of Theorem 7.1.}
    We begin with the existence.
By dilation we may assume that $|\partial\Omega|=1$. Let
 $f\in W^{1,p}(\partial\Omega)$.  Choose $f_k\in
    C^{\infty}_0(\rn{d})$ such that $f_k\rightarrow f$ in
    $W^{1,p}(\partial\Omega)$.  Let $u_k$ be the solution of the $L^2$
    regularity problem with data $f_k$.  Then $u_j-u_k$ is the
    solution of the $L^2$ regularity problem with data $f_j-f_k$
     Thus, using estimate (\ref{mainestimate}),
    we obtain
    \begin{eqnarray}
        \|(u_j)^*\|_{p}+\|(\nabla u_j)^*\|_{p} &\leq&
        C\|f_j\|_{W^{1,p}(\partial\Omega)}, \\
        \|(u_j-u_k)^*\|_{p}+ \|(\nabla u_j-\nabla u_k)^*\|_{p}
        &\leq& C\|f_j-f_k\|_{W^{1,p}(\partial\Omega)}. \label{diffRegEst}
    \end{eqnarray}
    It follows from estimate (\ref{diffRegEst}) that $u_j$ converges
 uniformly on any compact subset of $\Omega$. By interior estimates, this implies that
$u_j \to u$ and
$\nabla u_j \to \nabla u$ uniformly on any compact subset of $\Omega$ and
    $\mathcal{L}(u)=0$ in $\Omega$.
    By letting $j\rightarrow \infty$, as in the proof of Theorem 6.1, we obtain
    \begin{eqnarray}
\| (u)^*\|_p +\|( \nabla u)^*\|_{p}  &\leq&
        C\|f\|_{W^{1,p}(\partial\Omega)},  \\
 \| (u_k-u)^*\|_p + \|(\nabla u_k -\nabla u)^*
\|_{p} &\leq& C\|f_k-f\|_{W^{1,p}(\partial\Omega)} \label{7.5}.
    \end{eqnarray}
We point out that estimate (\ref{7.5}) implies that $u=f$ on $\partial\Omega$ in
the sense of nontangential convergence.
This follows from the same argument used in the
    proof of Theorem \ref{regImpliesDirichlet} for the existence of
nontangential limits.

To demonstrate the uniqueness, we fix $x\in \Omega$ and suppose that
$$
\left\{\begin{array}{ll}
    \mathcal{L}(u)=0 & \mbox{ in } \Omega, \\ u=0 & \mbox{ on }
    \partial\Omega, \\ (\nabla u)^* \in
    L^p(\partial\Omega).&\end{array}\right.
$$
We need to show
    that $u\equiv 0$ in $\Omega$.  To do this, let
    $G^x=G(x,y)=\Gamma (x-y)-w^x(y)$ in $\Omega$, where $w^x$ is the unique
solution of the $L^{p^\prime}$ Dirichlet problem in $\Omega$ with data
$\Gamma (x-y)$; i.e.,
    \begin{equation}
        \left\{\begin{array}{ll} \mathcal{L}(w^x)=0 & \mbox{ in }
        \Omega,
        \\ w^x=\Gamma^x & \mbox{ on } \partial\Omega, \\ (w^x)^*\in
        L^{p^\prime}(\partial\Omega).\end{array}\right.
    \end{equation}
Note that since $(w^x)^*\in L^{p^\prime}(\partial\Omega)$, we have
 $(G^x)^{*,\varepsilon}\in L^{p^\prime}(\partial\Omega)$
if $2\varepsilon<\text{dist}(x,\partial\Omega)$.

We now proceed as in the proof of Theorem 6.1.
It follows from (\ref{green}) that
    \begin{eqnarray}
        |u(x)| &\leq& \frac{C}{\varepsilon}
\int_{E_ \varepsilon} |G^x||\nabla u|\,dy + \frac{C}{\varepsilon^2}
        \int_{E_\varepsilon} |G^x||u| \,dy \nonumber \\
        &\leq& C\int_{\partial\Omega} (G^x)^{*,\varepsilon}\cdot(\nabla
        u)^*\,d\sigma \nonumber\\
        &\leq& C \|(G^x)^{*,\varepsilon}\|_{p^\prime} \|(\nabla u)^*\|_p,
        \label{uniqueEstimate}
    \end{eqnarray}
where we have used the observation that $|u(y)|\le C\varepsilon (\nabla u)^* (P)$ if
$y\in \gamma(P)\cap E_\varepsilon$.
Note that $\| (G^x)^{*, \varepsilon}\|_{p^\prime}\to 0$ as $\varepsilon\to 0$.
This follows easily from the dominated convergence theorem, since
    $(G^x)^{*,\varepsilon}\rightarrow 0$ a.e. as $\varepsilon\rightarrow
    0$ and $(G^x)^{*,\varepsilon} \leq (G^x)^{*,\varepsilon_0} \in
    L^{p^\prime}(\partial\Omega)$ for $\varepsilon \leq \varepsilon_0$.
    Thus we may conclude that $u\equiv 0$ in $\Omega$. This completes the proof.
\qed
\medskip

The rest of this section is devoted to the proof of Lemma \ref{mainLemma}.

\begin{lemma} \label{secLemma}
    Let $\Omega$ be a bounded Lipschitz domain with $|\partial\Omega|=1$ and $1<p<\infty$.
    Suppose that the $L^{p^\prime}$ Dirichlet problem
for $\mathcal{L}$ in $\Omega$ is uniquely solvable.
  Also suppose that $u\in
    C^\infty(\Omega)$,
$\nabla u$ exists a.e. on $\partial\Omega$,
and $u=0$ on $\Omega \backslash B(P,r)$ for some
$P\in \partial\Omega$. We further assume that for some $C_0>C_1>10$,
$\Omega\subset B(P,C_0r)$,
    $$
B(P,C_1r)\cap \Omega = B(P,C_1 r) \cap \big\{(x',x_d): x_d>\eta(x')\big\},
$$
 $\mathcal{L}(u)\in L^{p^\prime}(\Omega)$ and
$\left(\pdydx{u}{x_d}\right)^*\in L^{p^\prime}(\partial\Omega)$.
    Then,
    \begin{eqnarray*}
        \int_{ \partial\Omega} |\nabla u|^{p^\prime}\,d\sigma
        \leq C\int_{\partial\Omega}
        \left|\pdydx{u}{x_d}\right|^{p^\prime}\,d\sigma + C \int_{\Omega}
        |\mathcal{L}(u)|^{p^\prime}\,dx + C\int_{K} |\nabla
        u|^{p^\prime}\,dx,
    \end{eqnarray*}
 where $K$ is a compact subset of $\Omega$ and $C$ is a constant which may depend on $r$ and $K$.
\end{lemma}

\begin{proof} We may assume that $P=0$.
    Let $\widetilde{u}=u-\Gamma * \mathcal{L}(u)$ in $\Omega$.  Then
    $\mathcal{L}(\widetilde{u})=0$ in $\Omega$ and
    \begin{equation} \label{initialEstSecLemma}
        \int_{ \partial\Omega} |\nabla u|^{p^\prime}\,d\sigma
        \leq C\int_{\partial\Omega} |\nabla
        \widetilde{u}|^{p^\prime}\,d\sigma + C\int_{\partial\Omega}
        |\nabla \Gamma * \mathcal{L}(u)|^{p^\prime}\,d\sigma.
    \end{equation}

    \noindent Next, we estimate each of the terms on the right hand
    side of (\ref{initialEstSecLemma}).  We begin with the second
    term.  Note that
    \begin{eqnarray*}
        \int_{\partial\Omega} |\nabla \Gamma *
        \mathcal{L}(u)|^{p^\prime}\,d\sigma &\leq& C\int_{
        \Omega} |\nabla \Gamma * \mathcal{L}(u)|^{p^\prime}\,dx
+ C
        \int_{\Omega} |\nabla^2\Gamma *
        \mathcal{L}(u)|^{p^\prime}\,dx \\
        &\leq& C\int_{\Omega} |\mathcal{L}(u)|^{p^\prime}\,dx,
    \end{eqnarray*}
where we've used the Calder\'on-Zygmund estimates on the
    second term and the fractional integral estimates on the first
    term.  To estimate the first term in the right hand side of
    (\ref{initialEstSecLemma}),
we observe that by the square function estimates (e.g. see \cite{dkpv:areaintegral}),
    \begin{equation}\label{square-function-estimate}
        \int_{\partial\Omega} |\nabla
        \widetilde{u}|^{p^\prime}\,d\sigma \leq
 \int_{\partial\Omega} |(\nabla
        \widetilde{u})^*|^{p^\prime}\,d\sigma
        \leq C\int_{\partial\Omega} |S(\nabla
        \widetilde{u})|^{p^\prime}\,d\sigma + C\sup_{{K_1}} |\nabla
        \widetilde{u}|^{p^\prime},
    \end{equation}
where ${K_1}$ is a compact subset of $\Omega$ and $S(w)$ denotes the usual square
function of $w$, defined by using a regular family of nontangential cones.  Also note that
by interior estimates,
    \begin{eqnarray*}
        \sup_{{K_1}} |\nabla \widetilde{u}|^{p^\prime} &\leq&
        \int_{{K_2}} |\nabla \widetilde{u}|^{p^\prime}\,dx
        \leq \int_{{K_2}} |\nabla u|^{p^\prime}\,dx +
        C\int_{\Omega} |\mathcal{L}(u)|^{p^\prime}\,dx,
    \end{eqnarray*}
where $K_2\supset K_1$ is a compact subset of $\Omega$.

    It remains to estimate the term involving the
    square function in (\ref{square-function-estimate}). The key observation here is that
\begin{eqnarray}
S(\nabla\widetilde{u})& = &S(\nabla\Gamma*\mathcal{L}(u))
\le C\int_\Omega |\mathcal{L}(u)|\, dx \quad \text{ on }
\partial\Omega\setminus B(0,2r),\label{key-observation-1}\\
S(\nabla\widetilde{u}) &\le & C\, \widetilde{S}\left(\frac{\partial \widetilde{u}}{\partial x_d}\right)
+C \sup_{K_3} |\nabla \widetilde{u}|
\quad \text{ on } B(0,2r)\cap \partial\Omega,\label{key-observation-2}
\end{eqnarray}
where $K_3 \supset K_2$ is a compact subset of $\Omega$ and $\widetilde{S}(w)$ denotes the square
function of $w$, defined by using a regular family of nontangential cones which are slightly
larger than ones used for $S(w)$.
Estimate (\ref{key-observation-1}) follows easily from the assumption that
$u=0$ in $\Omega\setminus B(0,r)$.
The proof of (\ref{key-observation-2}) in the case of
harmonic functions on upper-half spaces may be found in \cite{stein:singularintegrals}.
It extends easily to the case of general second order elliptic systems in Lipschitz domains.

With estimates (\ref{key-observation-1})-(\ref{key-observation-2}) at our disposal, we have
    \allowdisplaybreaks{
    \begin{eqnarray}
         \lefteqn{\hspace{-0.25in}\int_{\partial\Omega}  |S(\nabla \widetilde{u})|^{p^\prime}\,d\sigma
        \leq \int_{B(0,2r)\cap \partial\Omega} |S(\nabla
        \widetilde{u})|^{p^\prime}\,d\sigma + \int_{\partial\Omega\setminus B(0,2r)}
        |S(\nabla \widetilde{u})|^{p^\prime}\,d\sigma}  \nonumber \\
        &\ & \leq C\int_{B(0,2r)\cap
        \partial\Omega} \left|
       \widetilde{S}\left(\pdydx{\widetilde{u}}{x_d}\right)\right|^{p^\prime}\,d\sigma +
      C\int_{\Omega} |\mathcal{L}(u)|^{p^\prime}\,dx
+  C \sup_{K_3} |\nabla \widetilde{u}|^{p^\prime} \nonumber \\
        &\ & \leq C\int_{\partial\Omega}
        \left|\left(\pdydx{\widetilde{u}}{x_d}\right)^*\right|^{p^\prime}\,d\sigma
        + C\int_{\Omega} |\mathcal{L}(u)|^{p^\prime}\,dx
+ C \sup_{K_3} |\nabla \widetilde{u}|^{p^\prime} \nonumber \\
        &\ & \leq C\int_{\partial\Omega} \left|
        \pdydx{\widetilde{u}}{x_d}\right|^{p^\prime}\,d\sigma + C\int_{\Omega}
        |\mathcal{L}(u)|^{p^\prime}\,dx + C \int_{{K}} |\nabla u|^{p^\prime}\,dx
        \label{LpSolvSecLem}
        \\
&\ & \leq C \int_{ \partial\Omega}
        \left|\pdydx{u}{x_d}\right|^{p^\prime}\,d\sigma +
        C\int_{\partial\Omega} |\nabla\Gamma *
        \mathcal{L}(u)|^{p^\prime}\,d\sigma \nonumber\\ && \qquad\qquad+ C\int_{\Omega}
        |\mathcal{L}(u)|^{p^\prime}\,dx \nonumber
  +C\int_{{K}} |\nabla u|^{p^\prime}\,dx \nonumber \\
        &\ & \leq C\int_{\partial\Omega}
        \left|\pdydx{u}{x_d}\right|^{p^\prime}\,d\sigma + C\int_{\Omega}
        |\mathcal{L}(u)|^{p^\prime}\,dx
+C \int_{{K}} |\nabla u|^{p^\prime}\,dx, \nonumber
    \end{eqnarray}
    }
where $K\supset K_3$ is a compact subset of $\Omega$.

     We point out that it is in estimate (\ref{LpSolvSecLem}) where
we used the solvability of the $L^{p^\prime}$ Dirichlet problem in $\Omega$. This is possible since
$\mathcal{L}(\frac{\partial \widetilde{u}}{\partial x_d}) =0$ in $\Omega$ and
    \begin{equation}\label{nontangential-estimate}
\left(\pdydx{\widetilde{u}}{x_d}\right)^*\in L^{p^\prime}(\partial\Omega).
\end{equation}
    To see (\ref{nontangential-estimate}), we only need to show that the radial maximal function
    \begin{equation}\label{radial-estimate}
\mathcal{M}\left(\pdydx{\widetilde{u}}{x_d}\right)(P)=\sup_{0<t<c}
    \left|\pdydx{\widetilde{u}}{x_d}(P+te_d)\right| \in
    L^{p^\prime}(B(0,2r)\cap\partial\Omega).
\end{equation}
Since
    $$
\mathcal{M}\left(\pdydx{\widetilde{u}}{x_d}\right) \leq
    \mathcal{M}\left(\pdydx{u}{x_d}\right) + \mathcal{M}(\nabla
    \Gamma * \mathcal{L}(u))
$$
and
    $\mathcal{M}\left(\pdydx{u}{x_d}\right) \leq
    \left(\pdydx{u}{x_d}\right)^*\in L^{p^\prime}(\partial\Omega)$,
estimate (\ref{nontangential-estimate}) follows from
    \begin{eqnarray}\label{radial-1}
        \int_{B(0,2r)\cap\partial\Omega} |\mathcal{M}(\nabla\Gamma *
        \mathcal{L}(u))|^{p^\prime}\,d\sigma
        \leq C\int_{\Omega} |\mathcal{L}(u)|^{p^\prime}\,dx.
    \end{eqnarray}
Finally we remark that the desired estimate (\ref{radial-1}) is a consequence of the inequality
$$
\int_{B(0,2r)\cap\partial\Omega} |\mathcal{M}(w)|^{p^\prime}\, d\sigma
\le C\int_\Omega |\nabla w|^{p^\prime} dx
+C\int_\Omega |w|^{p^\prime} dx
$$
for any $w\in C^1(\Omega)$. This complets the proof of Lemma \ref{secLemma}.
\end{proof}

We are now in a position to give

\noindent{\bf Proof of Lemma \ref{mainLemma}.} We may assume that $|\partial\Omega|=1$.
    Suppose that $f\in C_0^\infty(\rn{d})$ and
$$\left\{\begin{array}{ll} \mathcal{L}(u)=0 &
    \mbox{ in } \Omega, \\ u=f & \mbox{ on }
    \partial\Omega ,\\ (\nabla u)^*\in L^2(\partial\Omega). &
    \end{array}\right.
$$
By Theorem 2.1, $\nabla u$ has nontangential limits a.e. on $\partial\Omega$.
We will show that
\begin{equation}\label{conormal-estimate}
\left\|
    \pdydx{u}{\nu}\right\|_{p} \leq C \big\{ \|\nabla_t f\|_p +\| f\|_p\big\}.
\end{equation}
The nontangential maximal function estimate (7.1) follows
from (\ref{conormal-estimate}) and Theorem 2.1 as well as the Green's representation formula.
  To establish (\ref{conormal-estimate}), we first note that it suffices to
    consider the case $\mbox{supp}(f)\subset B(P_0,r)$, where $P_0\in \partial\Omega$ and
$B(P_0,C_1r)\cap \partial\Omega$ is given by the graph of a Lipschitz function after a possible
rotation.
  For    otherwise write $f=\sum_{j=1}^M f\varphi_j$ where $\varphi_j\in
    C^{\infty}_0(\rn{d})$ and $\sum \varphi_j =1$ on
    $\partial\Omega$.  Then, $u=\sum u_j$ where $u_j$ is the solution
    with data $f\varphi_j$ and we have
    \begin{eqnarray*}
        \left\|\pdydx{u}{\nu}\right\|_{p} \leq \sum_j \left\|
        \pdydx{u_j}{\nu}\right\|_{p}
        \leq C\sum_j \|f\varphi_j\|_{W^{1,p}(\partial\Omega)}
        \leq C\|f\|_{W^{1,p}(\partial\Omega)}.
    \end{eqnarray*}

    Let $g \in C^{\infty}_0(\rn{d})$ and
$$
\left\{ \begin{array}{ll}
    \mathcal{L}(w)=0 & \mbox{ in } \Omega, \\ w=g & \mbox{ on }
    \partial\Omega, \\ (\nabla w)^*\in L^2(\partial\Omega). &
    \end{array}\right.
$$
Using integration by parts we obtain
    $$\int_{\partial\Omega} \pdydx{u}{\nu}\cdot g\,d\sigma =
    \int_{\partial\Omega} \pdydx{u}{\nu}\cdot w\,d\sigma =
    \int_{\partial\Omega} u\cdot \pdydx{w}{\nu}\,d\sigma =
    \int_{\partial\Omega} f\cdot \pdydx{w}{\nu}\,d\sigma.$$
 By duality it suffices to show that
\begin{equation}\label{duality-estimate}
\left|\int_{\partial\Omega}
    f\cdot \pdydx{w}{\nu} \,d\sigma \right| \leq C\, \|f\|_{W^{1,p}(\partial\Omega)} \|g\|_{p^\prime}.
\end{equation}

    To this end we assume that $P_0=0$ and that
    $$
B(0,C_1 r)\cap \Omega = B(0,C_1r)\cap \big\{ (x',x_d): \ x_d>
    \eta(x')\big\},
$$
where $\eta:\rn{d-1}\rightarrow \rn{}$ is a
    Lipschitz function.  Next, choose $\psi \in
    C^{\infty}_0(B(0,4r))$ such that $\psi=1$ on
    $B(0,3r)$.  Now
    define
$$
\widetilde{w}(x',x_d) = -\int_{x_d}^{\infty}
    \psi(x',t)w(x',t)\,dt.
$$
Then $\pdydx{\widetilde{w}}{x_d} = \psi w$ in $\Omega$ and in particular,
 $\pdydx{\widetilde{w}}{x_d}=w$ in $B(0,3r)\cap
    \Omega$.
Also note that
    \begin{eqnarray*}
        (\mathcal{L}(\widetilde{w}))^\alpha &=& - \int_{x_d}^{\infty}
        (\mathcal{L}( \psi w))^\alpha\,dt \\
        &=& -\int_{x_d}^{\infty} a_{ij}^{\alpha\beta}\left\{
w^\beta D_iD_j\psi +D_iw^\beta D_j \psi
+D_j w^\beta D_i \psi \right\} \,dt.
    \end{eqnarray*}
It follows that on $B(0,2r)\cap \Omega$,
    \begin{eqnarray*}
        |\mathcal{L}(\widetilde{w})| \le C\sup_K |w|,
    \end{eqnarray*}
where $K\subset\subset \Omega$ is compact.  We now observe
    that
    \begin{eqnarray*}
        \int_{\partial\Omega} f \cdot \pdydx{w}{\nu}\,d\sigma &=&
        \int_{\partial\Omega} f^\alpha a_{ij}^{\alpha\beta} n_iD_jw^\beta\,d\sigma = \int_{\partial\Omega}
        f^\alpha a_{ij}^{\alpha\beta}n_iD_jD_d \widetilde{w}^\beta\,d\sigma \\
        &=& \int_{\partial\Omega} f^\alpha a_{ij}^{\alpha\beta}
        (n_iD_d-n_d D_i)D_j\widetilde{w}^\beta\,d\sigma +
        \int_{\partial\Omega}
        f^\alpha n_d[\mathcal{L}(\widetilde{w})]^\alpha\,d\sigma \\
        &=& - \int_{\partial\Omega} (n_iD_d-n_dD_i)f^\alpha\cdot
        a_{ij}^{\alpha\beta}D_j\widetilde{w}^\beta\,d\sigma + \int_{\partial\Omega}
        f^\alpha n_d[\mathcal{L}(\widetilde{w})]^\alpha\,d\sigma.
    \end{eqnarray*}
This implies that
    \begin{eqnarray*}
        \left| \int_{\partial\Omega} f\cdot \pdydx{w}{\nu}\,d\sigma
        \right| &\leq& C\|\nabla_t f\|_{p} \left\{
        \int_{B(0,r)\cap \partial\Omega} |\nabla
        \widetilde{w}|^{p^\prime}\,d\sigma\right\}^{\frac{1}{p^\prime}} + C\|f\|_{p}
        \sup_K |w| \\
        &\leq& C\|f\|_{W^{1,p}(\partial\Omega)} \left\{ \left\{ \int_{B(0,r)\cap
        \partial\Omega} |\nabla \widetilde{w}|^{p^\prime}\,d\sigma
        \right\}^{\frac{1}{p^\prime}} + \sup_K |w| \right\}.
    \end{eqnarray*}

    Note that since the $L^{p^\prime}$ Dirichlet problem in $\Omega$ is solvable, we have
    that $$\sup_K |w| \leq C\|(w)^*\|_{p^\prime} \leq C\|g\|_{p^\prime}.$$  Thus, we
    only need to show that
\begin{equation}\label{final-estimate}
\int_{B(0,r)\cap \partial\Omega}
    |\nabla \widetilde{w}|^{p^\prime}\,d\sigma \leq C\|g\|_{p^\prime}^{p^\prime}.
\end{equation}
Note that
    $\tilde{w}=0$ in $\Omega \backslash B(0,4r)$.  This allows us to apply
    Lemma \ref{secLemma} to obtain
    \begin{eqnarray}
        \int_{B(0,r)\cap \partial\Omega} |\nabla
        \widetilde{w}|^{p^\prime}\,d\sigma
 &\leq& C\int_{B(0,4r)\cap \partial\Omega}
        |w|^{p^\prime}\,d\sigma \nonumber \\ && \qquad+ C\int_{\Omega} |\mathcal{L}(\widetilde{w})|^{p^\prime}\,dx
        + C\int_{K} |\nabla \widetilde{w}|^{p^\prime}\,dx.\label{secLemmaMainLemma}
    \end{eqnarray}

To handle the second term  of (\ref{secLemmaMainLemma}), we use Hardy's inequality (see e.g. \cite{stein:singularintegrals}, p.272) to obtain
    \begin{eqnarray*}
\int_\Omega |\mathcal{L}(\widetilde{w})|^{p^\prime} dx &=&
        \lefteqn{\int_{B(0,4r)\cap \Omega} \left| \int_{x_d}^{\infty}
         \mathcal{L}(\psi w)
\,dt \right|^{p^\prime}\,dx}\\
 &\leq& C\int_{B(0,4r)\cap
        \Omega} |\mathcal{L}(\psi w)|^{p^\prime} \big\{ \text{dist}(x, \partial\Omega)\big\}^{p^\prime}\,dx \\
        &\leq& C\int_{\Omega} \big\{|\nabla w|^{p^\prime} +
        |w|^{p^\prime}\big\}\big\{ \text{dist}(x, \partial\Omega)\big\}^{p^\prime}\,dx \\
        &\leq& C\int_{\partial\Omega} |(w)^*|^{p^\prime}\,d\sigma \\
        &\leq& C\int_{\partial\Omega} |g|^{p^\prime}\,d\sigma.
    \end{eqnarray*}
By interior estimates, the last term in (\ref{secLemmaMainLemma}) is bounded by
    \begin{eqnarray*}
        C\sup_{\widetilde{K}} |w|^{p^\prime}
        \leq C\int_{\partial\Omega} |g|^{p^\prime}\,d\sigma.
    \end{eqnarray*}
where $\widetilde{K}\supset K$ is a compact subset of $\Omega$.  With this last estimate we complete the proof of estimate (\ref{final-estimate})
and hence the proof of Lemma \ref{secLemma}.
\qed

\bibliography{ks31}

\small
\noindent\textsc{Department of Mathematics,
University of Kentucky, Lexington, KY 40506}\\
\emph{E-mail address}: \texttt{jkilty@ms.uky.edu} \\

\noindent\textsc{Department of Mathematics,
University of Kentucky, Lexington, KY 40506}\\
\emph{E-mail address}: \texttt{zshen2@email.uky.edu} \\

\noindent \today

\end{document}